\documentclass[]{interact}
\usepackage{amsfonts,amsmath,amstext,amssymb,amsopn}
\usepackage{graphicx}
\usepackage{color}
\usepackage{url}
\usepackage{algorithm}
\usepackage{algorithmic}
\usepackage{natbib}
\bibpunct[, ]{(}{)}{;}{a}{}{,}

\usepackage{tabularx}
\usepackage{tabulary,booktabs, threeparttable, stackengine,array}
\usepackage{mathtools}

\usepackage[labelfont=bf]{subcaption}
\usepackage{bbold}
\usepackage{multirow}

\usepackage{mathtools}
\usepackage{cuted}

\usepackage[justification=centering]{caption}

\title{A Multi-stage Stochastic Programming Model for Adaptive Biomass Processing Operation under Uncertainty}

\author{
\name{B. Gulcan\textsuperscript{a}\thanks{CONTACT B. Gulcan. Email: bgulcan@g.clemson.edu}, Y. Song\textsuperscript{a}\thanks{CONTACT Y. Song. Email: yongjis@g.clemson.edu}, and S.~D. Eksioglu\textsuperscript{b}}
\affil{\textsuperscript{a}Department of Industrial Engineering, Clemson University, Clemson, SC 29634 USA; \textsuperscript{b}Department of Industrial Engineering, University of Arkansas, Fayetteville, AR 72701 USA}
}

\begin{document}
\newcommand{\exclude}[1]{}

\maketitle

\begin{abstract}
Variations of physical and chemical characteristics of biomass reduce equipment utilization and increase operational costs of biomass processing. Biomass processing facilities use sensors to monitor the changes in biomass characteristics. Integrating sensory data into the operational decisions in biomass processing will increase its flexibility to the changing biomass conditions. In this paper, we propose a multi-stage stochastic programming model that minimizes the expected operational costs by identifying the initial inventory level and creating an operational decision policy for equipment speed settings. These policies take the sensory information data and the current biomass inventory level as inputs to dynamically adjust inventory levels and equipment settings according to the changes in the biomass' characteristics. We ensure that a prescribed target reactor utilization is consistently achieved by penalizing the violation of the target reactor feeding rate. A case study is developed using real-world data collected at Idaho National Laboratory's biomass processing facility. We show the value of multi-stage stochastic programming from an extensive computational experiment. Our sensitivity analysis indicates that updating the infeed rate of the system, the processing speed of equipment, and bale sequencing based on the moisture level of biomass improves the processing rate of the reactor and reduces operating costs.
\end{abstract}


\section{Introduction}
The future of manufacturing lies in creating a connected digital system for physical assets (e.g., vehicles, equipment, machines, robots), people, and information systems. Integration of these different elements will result in smart manufacturing facilities and contribute to transitioning into Industry 4.0 with higher productivity, better quality control, and more flexible operations without the interference of a human operator \citep{FRANK201915}. One of the important elements in Industry 4.0 is sensing technology \citep{DALENOGARE2018383}. The use of sensory data allows the integration of sequentially revealed information into the decision-making process and increases the facilities' adaptability to dynamic system environments such as the uncertainty of characteristics of materials being processed. In this paper, we focus on biomass processing systems in biorefineries. Biorefineries use monitoring and sensing technologies to measure biomass characteristics before they are processed. Biorefineries typically use sensors to measure moisture level and weight. For example, Snetterton and Sleaford Renewable Energy Plants are two of the several commercial-size biomass processing plants in the UK~\citeyearpar{snettertonbiomass}. Both plants use automated bale handling systems that utilize sensors to collect data about moisture level of incoming bales. The data is used to control their operations.

Bales vary based on the physical/chemical characteristics of biomass (such as ash, moisture, carbohydrate contents, etc.). These variations are due to spatial variations in weather and soil, harvesting equipment used, etc. For example, consider the scenario when a number of bales of different types of feedstocks with different moisture levels and ash content are processed on the same equipment (e.g., a grinder). The resulting distribution of particle size and particle uniformity of processed biomass vary from one bale to the next. These variations lead to uneven biomass flow in a biomass processing system that could affect equipment utilization rate and lead to inconsistent conversion rates. The system may revise equipment processing speed to mitigate these variations. The sensory data enables these revisions to be adaptive to the changes in biomass' characteristics. Some examples of equipment processing speed are the rotational speed for the conveyors and discharge rate for the storage equipment.   The discharge rate for the storage equipment controls the inventory level inside it, which can help `smooth out' the biomass flow inside the biomass feeding system. This helps to achieve the prescribed system utilization target. The objective of this paper is to investigate the value of using sensory data and a stochastic optimization model to provide adaptive system control with uncertain biomass characteristics.

We propose a multi-stage stochastic programming (SP) model that minimizes the expected operational cost by identifying an optimal operational policy for equipment speed setting and infeed rate in a biomass processing system. Multi-stage SP is a well-known decision-making framework where the uncertainty (e.g., biomass characteristics) is revealed over time, and decisions are made sequentially based on the information revealed in each stage. Multi-stage SP models extend the two-stage SP models that provide optimal static strategic and tactical decisions by finding optimal adaptive decision policy in each time stage based on the uncertainty realized so far~\citep{Birge_1985}. The policy obtained from the multi-stage SP model will take sensory data and the current biomass inventory level as inputs and provide dynamic control on inventory levels by setting equipment processing speed and infeed rate in an adaptive fashion. 

We expect that the results of this study will facilitate the creation of an automated process control inside biorefineries that will increase productivity, result in better quality control and seamless and flexible operations. In addition, the results of this study will help biorefineries to develop strategies that increase reactor utilization and minimize costs; and to identify an optimal operational condition in the face of stochastic biomass characteristics. These outcomes will facilitate the commercial-scale generation of biofuels at a competitive cost. In the long run, these outcomes will strengthen the sustainable bioeconomy of the US,  enhance the security of energy supplies, reduce dependencies on fossil fuels, and reduce greenhouse gas (GHG) emissions. A strong and sustainable bioeconomy has additional socio-economic benefits, such as generating new green jobs, growth of rural economy and social stability~\citep{DOMAC2005Bioenergy,Sims2003Bioenergy,You2012Biofuel}.

The remainder of the paper is organized as follows. We conduct a brief literature review in Section~\ref{sec:LitReview}. In Section~\ref{sec:Problem}, we describe the  
problem setting and the modeling framework. In Section~\ref{sec:Solution}, we present our solution approach and implementation details. We conduct a case study and summarize numerical results and sensitivity analysis in Section~\ref{sec:Results}. We close with some concluding remarks in Section~\ref{sec:Conclusion}.

\section{Literature Review}\label{sec:LitReview}
Our study is closely related to two main streams of literature, namely multi-stage stochastic programming and optimization of biorefinery operations. 

\subsection{Multi-stage Stochastic Programming}\label{sec:Multi-stageLit}
Multi-stage SP is a well-known decision-making framework where the uncertain information is revealed over time, and the decisions are made sequentially based on the data available at each stage. Real-world applications of multi-stage SP models are rich in many areas such as energy~\citep{kharoufeh2018_ESS_in_MG, Multistage-hydrothermal2017, BRUNO2016979, siddig2019adaptive}, finance~\citep{Carino1994Multistage_Insurance, steinbach1999recursive, gulpinar2002multistage, JitkaPortfolia2009}, facility capacity planning ~\citep{Shabbir_SIP2003, CHEN2002781, GUPTA2014180, Singh2009Dantzig-WolfeDecomposition}, and transportation~\citep{ALONSO200047, Herer2006Transshipment, moller2008airline}.
While the main advantage of multi-stage models is to develop adaptive decision policies, it is computationally challenging to solve these models. The challenge is due to the uncertainty in the data and the nested structure of the multi-stage decision-making problem. A typical approach to handle the uncertainty is to approximate the underlying stochastic process using scenario trees~\citep{BL97}. However, as the number of decision stages in the planning horizon increases, the size of the scenario tree grows exponentially due to the nested structure~\citep{Shapiro2005ComplexitySP}. Many decomposition algorithms have been proposed in the literature to address this burden. \citep{DecompositionMultistage_Birge1985} extends the known two-stage L-shaped or Benders decomposition to multi-stage linear SP problems. The main idea of this nested Benders decomposition algorithm is to develop an outer approximation of the expected cost of the future stages (i.e., cost-to-go functions) by generating Benders cuts. Note that, as the number of stages increases, the nested Benders decomposition may become inefficient to solve the problem. \citep{Pereira_Pinto_1991} proposed the Stochastic Dual Dynamic Programming (SDDP) algorithm for multi-stage problems where the underlying random vectors are stage-wise independent and the expected cost-to-go functions are convex. Stage-wise independence means the probability distribution of the random variable that is observed in stage $t$, $\xi_t$, is independent of the realization of $\xi_{t-1}$. The stage-wise independence property allows one to define only a single expected cost-to-go function for every stage. SDDP is an iterative algorithm that has shown to have finite convergence with probability one under mild conditions~\citep{Shapiro_2011}, which maintains a cutting plane approximation for the expected cost-to-go functions. In our paper, we apply the SDDP algorithm to solve the proposed multi-stage SP model. We share more details about the algorithm in Section~\ref{sec:Solution}.

The literature shows the value of multi-stage SP (VMS) in many application areas. The VMS measures the relative advantage of multi-stage solutions over their two-stage counterparts. \citep{huang2009value} provides analytical bounds for the VMS in capacity planning under uncertainty. Similarly, \citep{XIE2018130} shows the VMS for strategic expansion of biofuel supply chain under uncertainty. \citep{Mahmutogullari2019VMS} presents the VMS in risk-averse unit commitment under uncertainty, a daily problem that arises energy market participation. 
This paper contributes to the current multi-stage SP literature by demonstrates the VMS for the optimization of manufacturing and processing facilities.

\subsection{Optimization of Biorefinery Operations}\label{sec:optBio_Lit}
Biorefinery operations include functions such as biomass preprocessing, biomass handling, and biomass storage. Works that focus on evaluating the design parameters of the equipment used cover most of the literature related to biorefinery operations \citep{Crawford_2016_Flowability_CornStover,DAI2012716}. Other studies analyze the energy consumption of equipment as a function of its design parameters and biomass characteristics \citep{osti_1369631, osti_1133890, yancey2015size}. These studies are limited in scope since they do not capture the interactions between equipment and the impact of equipment on the system's performance. In our paper, we take a system's approach and use a network flow model to characterize the flow of biomass over time during the planning. 

A number of studies use optimization models for optimal process control. These models minimize the energy consumption of individual equipment. For example, in \citep{NUMBI20161653}, the authors propose a deterministic, non-linear optimization model to determine a crusher's optimal operating parameters.  The model considers the time-to-use electricity tariff to achieve additional cost savings. Work by \citep{ZHANG20101929} presents an optimal control model for a series of conveyor belts. The authors evaluate two different operational structures: one with fixed conveyor speeds and another with conveyor speeds that are modified over time. The comparison of these models demonstrates the benefits of having the flexibility to adjust the operating conditions of a conveyor. To summarize, optimal control models are complex non-linear programs and thus they are limited to model only some part of the system. In our paper, we impose simplifying assumptions to avoid complicated non-linear and non-convex constraints in the proposed multi-stage SP model. 

A few studies present models for optimal design of operations in a biorefinery. \citep{Viet2012BiorefOpt}, for example, proposes a two-step strategy to optimizing a biorefinery's design and operations. This model identifies an optimal biorefinery configuration for a given set of biomass feedstocks and available conversion technologies. The first step is to determine which intermediary chemicals can be produced using the available feedstocks. The second step is a network flow optimization problem that identifies the optimal combination of intermediary chemicals and corresponding conversion technologies that minimize the total production cost. Work by \citep{ZONDERVAN2011BiorefOpt} uses a mixed-integer and non-linear network optimization model to identify an optimal design for a biorefinery. The model identifies process sequences to optimize the production of a set of biofuels and bioproducts. Both works, \citep{Viet2012BiorefOpt} and \citep{ZONDERVAN2011BiorefOpt}, use deterministic models since they consider fixed biomass characteristics, flow rates, and yields.  Different from these works, our research models the uncertainty of biomass moisture level. In a recent work, we proposed a chance-constrained SP model in~\citep{Gulcan2021} that identifies operating conditions of equipment and inventory level to maintain a continuous flow of biomass to the reactor. This work considers the uncertainty of biomass moisture content,  particle size uncertainties, and random equipment failure. This chance-constrained SP model results in static policies that cannot incorporate newly learned information via sensory regarding of biomass characteristics. In this paper, our proposed multi-stage SP model will enable the integration of the sensory data into the decision-making process and increase the facilities' adaptability to the uncertain characteristics of biomass. 

\section{Problem Description and Mathematical Formulation}\label{sec:Problem}
\begin{figure}[h]
	\centering
	\includegraphics[width=\textwidth]{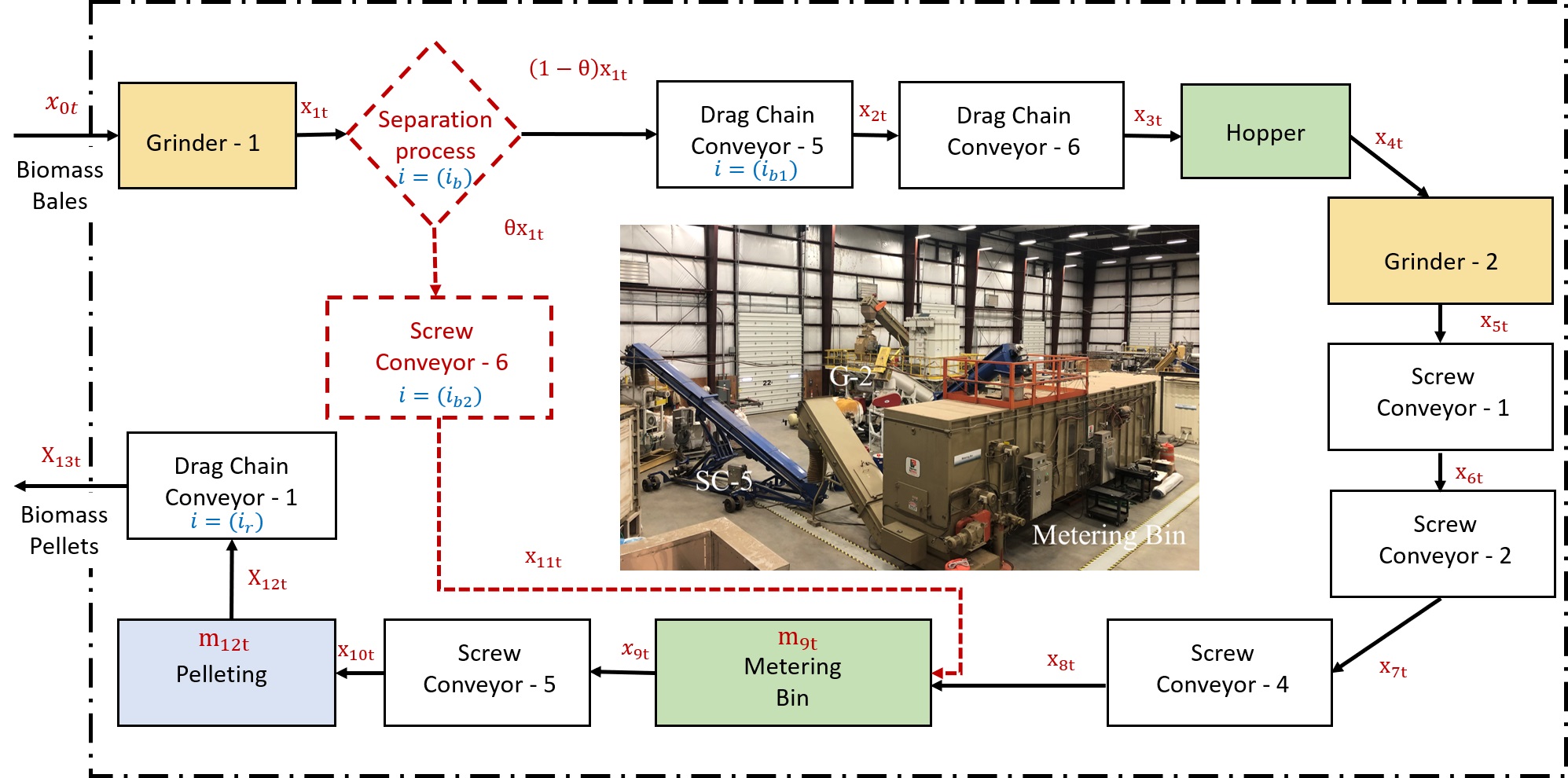}
	\caption{Biomass Feeding Processes of PDU}
	\label{fig:PDU}
\end{figure}

Figure \ref{fig:PDU} presents the flowchart of the different processes of the Processing Development Unit (PDU), a full-size, fully-integrated   feedstock preprocessing system at the Idaho National Laboratory (INL). In this figure, notations $x_{1t}, x_{2t}, ...$ represent the flow of biomass in stage $t.$ Bales of biomass enter the system, one at a time, via a conveyor belt. Biomass is processed in Grinder 1. Processed biomass undergoes a screening procedure that separates biomass based on particle size. Large biomass particles are transported via conveyors to Grinder 2 for further processing. Small biomass particles bypass Grinder 2 and a conveyor transports them to the Metering bin for storage. Biomass  processed in Grinder 2 also is stored in the Metering bin. Conveyors connected to the Metering bin transports the stored biomass to the pelleting machine. Pelleted biomass is fed to the reactor. The maximum infeed rate to the equipment and target feeding rate of the reactor are provided in `dry ton per hour' ($dt/hr$).

We group equipment in this system based on their tasks into processing, ${\bf E}^p$, transportation, ${\bf E}^r$, and storage, ${\bf E}^m$. A processing equipment converts biomass from its original format (e.g., baled, log, or coarse-shredded) to its final format (e.g., ground biomass and pellet). Grinders 1 and 2 and the pelleting mill are processing equipment. Most equipment of PDU are transportation equipment, such as conveyors. Storage equipment, such as the metering bin, store in-process inventory.

We use a network flow model to characterize the flow of biomass within the system over time during the planning horizon $\mathcal{T}$. Let $\bf G=(\bf N, \bf A)$, with a node set $\bf N$ and arc set $\bf A$, represent the network structure of the system. The set of nodes represents the equipment of the PDU (Figure \ref{fig:PDU}), and the set of arcs represents the flow of biomass from one equipment to the next. The following is a list of notations we use.

\begin{table}[H]
	\scriptsize
		\begin{tabular}{ll}
			\toprule
			\multicolumn{2}{c}{\textbf{SETS:}}\\
		    ${\bf M}$ & The set of moisture levels of biomass, 
	        ${\bf M} := \{LOW,\ MED,\ HIGH\}$.\\
		    $\mathcal{T}$ & The set of stages in the planning horizon, $\mathcal{T} := \{1,2,\ldots, T\}$\\
		    ${\bf E}^p$ & The set of processing equipment.\\
		    ${\bf E}^r$ & The set of transportation equipment.\\
		    ${\bf E}^m$ & The set of storage equipment.\\
		    ${\bf N} =  {\bf E}^p \cup {\bf E}^r \cup {\bf E}^m$\\
		    $\delta_i^+$ & The set of equipment that feed storage equipment $i \in E^m$.\\
		    $\delta_i^-$ & The set of equipment that are fed by storage equipment $i \in E^m$.\\
			\midrule
			\multicolumn{2}{c}{\textbf{PARAMETERS:}}\\
			$A$ & Node arc incidence matrix of $\bf G$.\\
			$A^{\prime}$ & Node arc incidence matrix of $\bf G^{\prime}(\bf {\bf E}^r, \bf A^{\prime})$.\\
		    $\mathbf{I}_0$ & Initial inventory level of the storage equipment (in dt).\\
			$\mathbf{c}^h$ & Inventory holding cost for storage equipment (in \$/dt).\\
	        $\mathbf{c}^p$ & Penalty cost of not satisfying the target utilization of the reactor (in \$h/dt).\\
	        $\kappa_t$ & Moisture level of the biomass bales in stage t ($\kappa_t \in M$).\\	
			$\bar{v}_t(\kappa_t)$ & Upper bound of equipment processing speed, with respect to the moisture level.\\ 
			$\bar{\iota}(\tilde{m}_t)$ & Inventory holding capacity of the storage equipment (in dt).\\
			$r$ & The target reactor feeding rate (in dt/h).\\
			$v_{1t}(\kappa_t)$ & The system feeding rate, with respect to the moisture level in stage $t$ (in m/hr).\\ 
			$i_r$ & Index of the last equipment that feeds the reactor\\
			\midrule
			\multicolumn{2}{c}{\textbf{RANDOM PARAMETERS:}}\\
			$\tilde{m}_t$ & Moisture content of the biomass bales in stage $t$. (in \%)\\
			\midrule
			\multicolumn{2}{c}{\textbf{DECISION VARIABLES:}}\\
			$\mathbf{V}_t := \{V_{it}\}_{i\in {\bf N}}$,  & Equipment processing speed (in centimeters or rotations per minute) in stage $t \in \mathcal{T}$.\\ 
			$\mathbf{I}_t := \{I_{it}\}_{i\in {\bf E}^m}$ & Inventory level in stage $t \in \mathcal{T}$ (in  dt).\\
			$\mathbf{X}_t := \{X_{it}\}_{i\in {\bf N}}$ & Biomass flow in stage $t$ (in dt).\\
			$p_t$ & Shortfall of biomass to achieve the target feeding of the reactor in stage $t \in \mathcal{T}$ (in dt).\\  
			\bottomrule
			\end{tabular}
\end{table}%

In a multi-stage SP setting, the uncertain information is revealed gradually over time, and the decision-maker may adapt their decisions accordingly. The inventory level $I_t$ in the storage equipment characterizes the state of the biomass processing system and links stages of the planning horizon. In the context of our problem, sensors measure the moisture content of the biomass, $\tilde{m}_t$, that will be processed over the processing time. Biomass moisture is one of the main factors that impact the biomass flow in the system by affecting biomass density~\citep{ZhouCornStover2008},~\citep{Crawford_2016_Flowability_CornStover}. 

We use biomass density to calculate the flow of biomass in the system. In addition, transportation and storage equipment have volumetric capacities. The changes in biomass density impact the amount of  biomass that a piece of equipment can transport or store. Thus, it is important that the model considers biomass density and how uncertain biomass moisture impacts the density.

Work by \citep{yancey2015size} shows that increase of moisture content results in a higher clogging rate. Thus, it is a common practice to lower  processing rates of equipment when processing biomass with high moisture level. As a result, the  moisture level of biomass affects both the system infeed rates, $\mathbf{V}_t$, and biomass flow, $\mathbf{X}_t$. Thus, making decisions by considering the information about the biomass moisture content can improve system's performance.

The proposed multi-stage SP model is created based on the following key definitions:

\begin{itemize}
    \item \emph{Stages}: We define a stage as the time period during which the system  processes one biomass bale. At each stage, the exact numerical values of moisture content, $\tilde{m}_{t}$, will be known via sensors' readings. 
    
    \item \emph{State variables:} The inventory level ($\mathbf{I}_t$) in the storage equipment characterizes the state of the biomass processing system, which is carried over from one stage to the next.
    
    \item \emph{Actions}: At each stage, we can change the processing speed of different equipment $\mathbf{V}_{t}$ according to the realization of biomass moisture level (provided by sensors' readings) and the current inventory levels.
\end{itemize}

The following assumptions are employed in our model:
\begin{itemize}
    \item The random variables $\tilde{m}_t$'s are stage-wise independent. This is a realistic assumption in that biomass is stored by bales -- the distribution of the random moisture content of the next bale is independent of the realized moisture level of the current bale.
    \item The system will process a single bale at each stage. This is an assumption we make to clearly model the structure of our problem and define a stage. In Section~\ref{sec:Sensitivity_Granularity}, we relax this assumption and consider different stage definitions. 
    \item The system infeed rate, i.e., the processing rate of the first equipment, $v_{1t}$, is  given as a problem parameter. If $v_{1t}$ is introduced as a decision variable, the problem becomes non-linear and non-convex due to the flow and inventory balance constraints in the storage equipment. The resulting non-convex multi-stage SP model will be computationally much more challenging to solve. In addition, the value of $v_{1t}$ does not impact the total amount of biomass flowing into the system in a stage but only affects the duration of a stage. 
    \item We only consider the uncertainty of moisture content of biomass bales, and we assume that given the moisture content, the biomass density value after being processed in grinders 1 and 2 is calculated via a deterministic function.

    \item We do not consider the blending of  biomass of different moisture content in the storage equipment. As a result, the density of  biomass in the storage unit is calculated using the realizations of moisture content of the current stage. Considering the biomass blending inside the storage equipment will lead to a non-convex multi-stage SP model that is very challenging to solve.
\end{itemize}

We are now ready to present the proposed multi-stage SP model. First, we introduce a generic nested formulation~\eqref{eqn:Nested}. Let $\mathbf{Y}_t\ := \ \{V_t, X_t, p_t\}$ be the vector of local stage variables (i.e., control variables) for each stage $t\in \mathcal{T}$, and let $\mathcal{Y}_t(\mathbf{I}_{t-1}, \tilde{m}_t)$ be the feasible region for state variables $\mathbf{I}_t$ and local variables $\mathbf{Y}_t$ given the previous inventory level $\mathbf{I}_{t-1}$ and a realization of random vector $\tilde{m}_t$, a nested formulation of the multi-stage SP model can be presented as:

\begin{align}\label{eqn:Nested}
    \min_{(\mathbf{I}_1, \mathbf{Y}_1)\in \mathcal{Y}_1(\mathbf{I}_0, \tilde{m}_1)}  \mathbf{c}^{h \top} \mathbf{I}_1 + \mathbf{c}^{p \top}\mathbf{Y}_1  + \mathbb{E} \left. \bigg\lbrack \right. & \min_{(\mathbf{I}_2,\mathbf{Y}_2) \in \mathcal{Y}_2(\mathbf{I}_1, \tilde{m}_2)} \mathbf{c}^{h \top} \mathbf{I}_2 +  \mathbf{c}^{p \top}\mathbf{Y}_2   \nonumber \\
    & + \left. \mathbb{E} \left \lbrack \cdots + \mathbb{E} \left \lbrack  \min_{(\mathbf{I}_T,\mathbf{Y}_T) \in \mathcal{Y}_T(\mathbf{I}_{T-1}, \tilde{m}_T)} \mathbf{c}^{h \top} \mathbf{I}_T +  \mathbf{c}^{p \top}\mathbf{Y}_T  \right \rbrack \right \rbrack \right \rbrack,
\end{align} 

where the initial realization of the moisture content, $\tilde{m}_1$, is assumed to be deterministic. 

Unfortunately,  formulation~\eqref{eqn:Nested} is computationally intractable due to the nested optimization posed by the sequential nature of the decision-making structure. This challenge can be addressed by using a dynamic programming reformulation. Under the aforementioned stage-wise independence assumption, dynamic programming formulation can be solved efficiently. 

Now, we formulate the problem to be solved in each stage $t \in \mathcal{T}$ using the dynamic programming approach. Given a realization of the random vector $\tilde{m}_t$ and the previous inventory level $\mathbf{I}_{t-1}$, the optimization problem to be solved in stage $t$ is given by:

\begin{subequations}\label{MSP-biorefinary}
\begin{align}
Q_t(\mathbf{I}_{t-1},\tilde{m}_t):= \min \ &\resizebox{0.24\textwidth}{!}{$ \mathbf{c}^{h \top} \mathbf{I}_t + c^pp_t  + \mathcal{Q}_{t+1}(\mathbf{I}_{t})$} \label{eqn:multi-stage-obj} \\
\text{s.t. } &   A^{\prime}\mathbf{X}_t = 0,  \label{eqn:multi-stage-Flow1}\\
& \mathbf{X}_t = g(\mathbf{V}_t, \tilde{m}_t), \label{eqn:multi-stage-Flow2} \\
&  \mathbf{X}_t \le h(\mathbf{V}_t, \tilde{m}_t), \label{eqn:multi-stage-Flow3} \\
& \mathbf{I}_{it} = \mathbf{I}_{i,t-1} + \sum_{j \in \delta_i^+} X_{j,t}\nonumber\\
&  - \sum_{l \in \delta_i^-} X_{l,t}, \ \forall i \in E^m, \label{eqn:multi-stage-Inv}\\
& p_t \geq r - X_{i_r,t}, \label{eqn:Shortfall1}\\
& p_t \geq 0, \label{eqn:Shortfall2}\\
& 0 \leq \mathbf{V}_t \leq \bar{v}_{\kappa_t},
\label{eqn:multi-stage-Speed_bound}\\
& 0 \leq \mathbf{I}_t \leq \bar{\iota}(\tilde{m}_t), \label{eqn:multi-stage-InvLimit}
\end{align}
\end{subequations}

where 

\begin{equation*}
    \mathcal{Q}_{t+1}(\mathbf{I}_{t}) := \mathbb{E}\left \lbrack Q_{t+1}(\mathbf{I}_{t},\tilde{m}_{t+1}) \right \rbrack, \ \forall t \neq T,
\end{equation*}
\begin{equation*}
  \text{and }  \mathcal{Q}_{T+1}(\mathbf{I}_{T}) = 0. 
\end{equation*}

In formulation~\eqref{MSP-biorefinary}, the objective is to minimize the expected total inventory holding cost and the expected penalty of not achieving the target reactor feeding rate $r$. The latter objective helps ensure a consistently high utilization of the reactor, as the reactor is the most expensive equipment in the biorefinery. Constraint \eqref{eqn:multi-stage-Flow1} represents the flow balance for transportation equipment. Constraint \eqref{eqn:multi-stage-Flow2} calculates  biomass flow from storage and processing equipment. This constraint is a function of processing speed, inventory level, and moisture content of the biomass. Constraint \eqref{eqn:multi-stage-Flow3} represents the upper limit on the amount of biomass flow from equipment $i\in {\bf N}$. Functions $g(\cdot)$ and $h(\cdot)$ are linear functions with respect to equipment geometry, equipment processing speed, and biomass density, which depend on the biomass moisture content $\tilde{m}_{t}$. Below we share examples of $g(\cdot)$ and $h(\cdot)$: 
\[ X_{it} = \gamma_i d_{it}(\tilde{m}_{t}) V_{it}.\]
\[ X_{it}\leq \gamma_i d_{it}(\tilde{m}_{t}) V_{it}.\]
Here, $\gamma_i$ represents the geometry of equipment $i$ and $d_{it}(\tilde{m}_{t})$ represents  biomass density as a function of moisture content.

 The inventory balance constraints \eqref{eqn:multi-stage-Inv} link the inventory levels in successive stages. There can be multiple equipment which feed and are fed by the storage equipment ($\delta_i^-$ and $\delta_i^+$, respectively). Constraints \eqref{eqn:Shortfall1} and \eqref{eqn:Shortfall2} calculate the shortfall of achieving the target reactor feeding rate, which is penalized in the objective. Constraints \eqref{eqn:multi-stage-Speed_bound} and \eqref{eqn:multi-stage-InvLimit} set bounds on the processing speed of equipment and amount of inventory stored, respectively.
 
 \subsection{Two-stage Approximations}\label{sec:2stage}
 
\begin{figure}[H]
	\centering
	\begin{subfigure}{0.48\textwidth}
	    \centering
	    \includegraphics[width=0.75\textwidth]{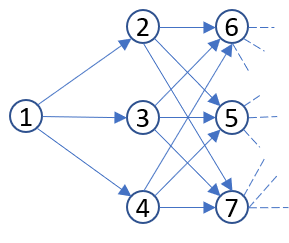}
	\caption{Multi-stage SP.}
	\label{fig:SC-tree}
	\end{subfigure}
	\hfill
	\begin{subfigure}{0.48\textwidth}
	    \centering
    	\includegraphics[width=\textwidth]{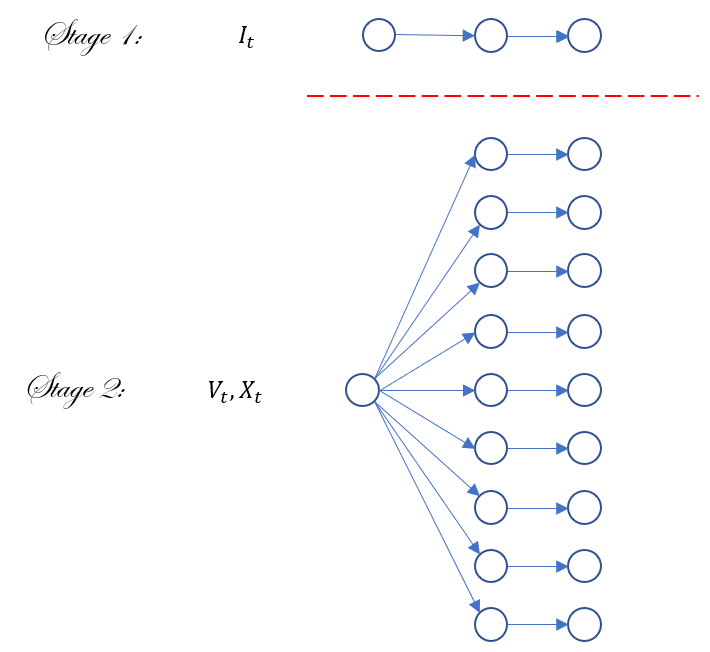}
    	\caption{Scenario tree illustration for two-stage approximation to multi-stage SP.}
    	\label{fig:SC-tree_2stage}
	\end{subfigure}
	\caption{Scenario tree illustrations}
\end{figure}

 The computational complexity of multi-stage SPs grows exponentially with the increase of the number of stages~\citep{Shapiro2005ComplexitySP}. Two-stage SP models are often used to approximate multi-stage SP models by making the state variables static over time as opposed to allowing them to be adaptive to dynamically revealed information. To illustrate this approximation, we use the scenario tree in Fig.~\ref{fig:SC-tree} to demonstrate the structure of the multi-stage SP model. Note that due to our assumption of stage-wise independent uncertainty, the scenario tree structure illustrated in the figure is referred to as a `recombining' scenario tree in the literature.

The following setting creates the two-stage approximation:
\begin{itemize}
\item In the first stage (prior to any realization of uncertainty), the inventory levels $\mathbf{I}_t$'s are determined for all $t\in \mathcal{T}$.
\item In the second stage, given the inventory levels for all $t\in \mathcal{T}$, the equipment processing speeds $\mathbf{V}_t$ and biomass flow $\mathbf{X}_t$ are determined based on the realization of uncertainty at each $t\in \mathcal{T}$. 
\end{itemize}

Figure \ref{fig:SC-tree_2stage} represents the scenario tree structure for the two-stage approximation. As we set $\{\mathbf{I}_t\}_{t\in \mathcal{T}}$ in the first stage, there is no distinction between different sample paths (a sample path is a sequence of realizations of random variables from $t = 1$ to $t = T$). Note that, $\{\mathbf{I}_t\}_{t \in \mathcal{T}}$ corresponds to the state variables that link different stages together. When $\{\mathbf{I}_t\}_{t\in \mathcal{T}}$ is given, there is no link between problems at different nodes in the second stage, and therefore, we can decompose the second stage by solving a separate problem at each node in the scenario tree given a first-stage solution. As we enforce $\{\mathbf{I}_t\}_{t\in \mathcal{T}}$ to be static decisions, i.e., they are enforced to be identical for all sample paths, this two-stage approximation is a \emph{restriction} of the multi-stage SP model. The key advantage of the proposed multi-stage SP model is to adapt decisions with respect to newly revealed information at each stage.

Given a set $S$ of sample paths (obtained via i.i.d samples from $\{\tilde{m}_{t}\}_{t\in \mathcal{T}}$), where $\tilde{m}^s_t$ gives the realization of $\tilde{m}_{t}$ on sample path $s$, the two-stage approximation model is presented below. Note that decision variables with a superscript are second-stage decision variables, while the ones without any superscript are first-stage decision variables. We will compare this approximation model with~\eqref{MSP-biorefinary} to assess the trade-off between the benefit of adaptive decision-making and the computational effort in our numerical experiments. 
\begin{subequations}\label{Two-stage_Approx}
\begin{align}
\min \ & \sum_{t \in \mathcal{T}} \mathbf{c}^{h \top} \mathbf{I}_t + \displaystyle \frac{1}{|S|}\sum_{s\in S} \sum_{t \in \mathcal{T}} c^p p^s_t  \\
\text{s.t. } & \mathbf{I}_t \ge 0, \ \forall t \in \mathcal{T} \\
&   A^{\prime}\mathbf{X}^s_t = 0, \ \forall t \in \mathcal{T}, \ \forall s\in S \\
& \mathbf{X}^s_t = g(\mathbf{V}^s_t, \tilde{m}^s_t), \ \forall t \in \mathcal{T}, \ \forall s\in S \\
&  \mathbf{X}^s_t \le h(\mathbf{V}^s_t, \tilde{m}^s_t),  \ \forall t \in \mathcal{T}, \ \forall s\in S \\
& \mathbf{I}_{it} = \mathbf{I}_{i,t-1} + \sum_{j \in \delta_i^+} X^s_{j,t} - \sum_{l \in \delta_i^-} X^s_{l,t}, \nonumber \\
&\hspace{.4in} \forall i \in E^m, \ \forall t \in \mathcal{T}, \ \forall s\in S \\
& p^s_t \geq r - X^s_{i_r,t}, \ \forall t \in \mathcal{T}, \ \forall s\in S \\
& p^s_t \geq 0, \ \forall t \in \mathcal{T}, \ \forall s\in S \\
& 0 \leq \mathbf{V}^s_t \leq \bar{v}_{\kappa_t},\ \forall t \in \mathcal{T}, \ \forall s\in S \\
& 0 \leq \mathbf{I}_t \leq \bar{\iota}(\tilde{m}^s_t), \ \forall t \in \mathcal{T},\ \forall s \in S.
\end{align}
\end{subequations}

 \subsection{Deterministic Model: Mean-value Problem}\label{sec:Deterministic}
 Solving stochastic programs can be computationally challenging. Many decision-makers prefer to solve simpler versions of the real-world problems. In this section, we introduce the \emph{mean-value (MV) problem}, which is deterministic and simple to use. We demonstrate the \emph{value of stochastic solutions(VSS)} by conducting an out-of-sample evaluation of the solutions to the MV problem and the solutions to stochastic model.
 
 \begin{subequations}\label{MVproblem}
\begin{align}
\min \ & \sum_{t \in \mathcal{T}}( \mathbf{c}^{h \top} \mathbf{I}_t +  c^p p_t ) \\
\text{s.t. } & \mathbf{I}_t \ge 0, \ \forall t \in \mathcal{T} \\
&   A^{\prime}\mathbf{X}_t = 0, \ \forall t \in \mathcal{T}\\
& \mathbf{X}_t = g(\mathbf{V}_t, \bar{m}_t), \ \forall t \in \mathcal{T} \\
&  \mathbf{X}_t \le h(\mathbf{V}_t, \bar{m}_t),  \ \forall t \in \mathcal{T} \\
& \mathbf{I}_{it} = \mathbf{I}_{i,t-1} + \sum_{j \in \delta_i^+} X_{j,t} - \sum_{l \in \delta_i^-} X_{l,t}, \nonumber \\
&\hspace{.4in} \forall i \in E^m, \ \forall t \in \mathcal{T} \\
& p_t \geq r - X_{i_r,t}, \ \forall t \in \mathcal{T} \\
& p_t \geq 0, \ \forall t \in \mathcal{T} \\
& 0 \leq \mathbf{V}_t \leq \bar{v}_{\kappa_t},\ \forall t \in \mathcal{T} \\
& 0 \leq \mathbf{I}_t \leq \bar{\iota}(\bar{m}_t), \ \forall t \in \mathcal{T}.
\end{align}
\end{subequations}

\section{Solution Methodology: Stochastic Dual Dynamic Programming}\label{sec:Solution}
In this section, we discuss the solution approaches for solving the proposed multi-stage SP model and implementation details. Recall that the SDDP algorithm proposed by~\citep{Pereira_Pinto_1991} is a popular approach for multi-stage SP with stage-wise independent uncertainty and convex expected cost-to-go functions. Both of these two assumptions are satisfied in model~\eqref{MSP-biorefinary} as it only involves linear constraints and continuous decision variables. The SDDP algorithm maintains and iteratively updates a Benders-type cutting-plane approximation for the expected cost-to-go function $Q_t(\cdot)$ until a termination criterion is met. There are two main steps in each iteration of the SDDP algorithm: a backward pass and a forward pass.

To illustrate the SDDP algorithm, let us use the following simplified notation for the stage $t$ problem:
\begin{align*}
Q_t(\mathbf{I}_{t-1},\xi_t) := \min \ &  \mathbf{c}^{\top} \mathbf{I}_t +  \mathcal{Q}_{t+1}(\mathbf{I}_{t})  \\
\text{s.t. } &   A_t  \mathbf{I}_{t} = b_t - B_t \mathbf{I}_{t-1},
\end{align*}

where $\xi_t := (c, A_t, b_t, B_t)$ represents the (random) problem data. In general, some or all of the data $(c, A_t, b_t, B_t)$ can be subject to uncertainty. In the context of our problem, $c$ is deterministic and $(A_t, b_t, B_t)$ are subject to uncertainty according to their dependence on $\tilde{m}_t$.

\paragraph{Forward pass} Sample $M << N$ paths, where $N$ is the total number of sample paths in the scenario tree. Based on the current approximation of the expected cost-to-go functions, $\hat{\mathcal{Q}}_{t+1}(\cdot)$, and the candidate solution for the state variable in the previous stage, $\mathbf{\bar{I}}_{t-1}$, solve the stage $t$ problem to obtain a candidate solution, $\mathbf{\bar{I}}_{t}$, for each stage $t \in \mathcal{T}$ on each of the $M$ sample paths. 

\paragraph{Backward pass} Let $N_t$ be the number of realizations of the random variable $\xi_t$ considered in the scenario tree at stage $t$. Let $\xi^j_{t} := (c,A^j_t, b^j_t, B^j_t)$ for $j \in \{1,\dots, N_t\}$ correspond to realization of the $j^{th}$ node in stage $t$. In a backward pass, the stage $t$ problem, shown below, is solved for all $N_t$ realizations of $\xi_t$:   

\begin{align*}
\underline{Q}_t(\mathbf{\bar{I}}_{t-1},\xi^j_t) := \min \ &  \mathbf{c}^{\top} \mathbf{I}_t  + \hat{\mathcal{Q}}_{t+1}(\mathbf{I}_{t})  \\
\text{s.t. } &   A^j_t  \mathbf{I}_{t} = b^j_t - B^j_t \mathbf{\bar{I}}_{t-1} \ \ \  (\pi^j_t),
\end{align*}
where $\pi^j_t$ is the associated optimal dual vector. After solving all $N_t$ subproblems, we have the following Benders cut:

\[
\theta_{t-1} \geq \frac{1}{N_t} \sum^{N_t}_{j = 1} \underline{Q}_t(\mathbf{\bar{I}}_{t-1},\xi^j_t) - \frac{1}{N_t} \sum^{N_t}_{j = 1} B^j_t \pi^j_t (\mathbf{I}_{t-1} - \mathbf{\bar{I}}_{t-1})
\]

This cut is added to the collection of cutting planes to improve the approximation $\hat{\mathcal{Q}}_t(\cdot)$. 

\paragraph{Termination condition} As each of the expected cost-to-go function approximation $\{\hat{\mathcal{Q}}_t(\cdot)\}_{t=2}^T$ is a lower approximation to the true function, at each iteration, the optimal value of the first-stage problem with $\hat{\mathcal{Q}}_2(\cdot)$ provides a deterministic lower bound. After each forward pass, we can calculate a statistical upper bound for the problem. For each sample path $j=1,\dots, M$, let $v_j$ be the sum of immediate cost associated with the trained optimal policy over all stages. Let $\bar{v} := \frac{1}{M}\sum^M_{j = 1} v_j$ be the sample average and $\sigma^2_v := \frac{1}{M-1}\sum^M_{j = 1} (v_j - \bar{v})^2$ be the sample variance. The statistical upper bound is calculated as follows:
\[
    v^{UB} := \bar{v} + \frac{z_\alpha \sigma_v}{\sqrt{M}}.
\]
Here $z_\alpha$ denotes the $(1 - \alpha)$ quantile of the standard normal distribution, and $v^{UB}$ gives an upper bound for the optimal value of the true problem with a confidence of $(1-\alpha)$. We used $\alpha = 0.025$, which corresponds to the $95\%$
confidence interval. One meaningful termination condition is checking the gap between $\bar{v}$ and $v^{UB}$ and comparing it to a prescribed threshold parameter~\citep{Shapiro_2011}. Another termination condition that might be computationally more efficient in practice is checking the progress of the improvement of deterministic lower bound over iterations. 
The algorithm can be terminated if this lower bound fails to improve by more than a user-specified threshold $\epsilon$ for more than a certain number of consecutive iterations $\sigma$. For example, let $\epsilon = 10^{-4}$ and $\sigma = 10$, if for the last $10$ iterations, the lower bound does not improve by more than $10^{-4}$, then the SDDP algorithm terminates.

\paragraph{SDDP Implementation} 

We used the Julia package \citep{dowson_sddp.jl} that implements the SDDP algorithm to solve the multi-stage SP model. We considered 500 realizations of $\tilde{m}_t$ in each stage $t$ to solve the multi-stage SP model via the SDDP algorithm. The package refers to a solution of the model as the trained policy. Various termination conditions are implemented in this Julia package. We chose to use the bound stalling with a tolerance of $10^{-4}$ and 50 previous iterations to check. Once a policy is obtained by solving the multi-stage SP model, we evaluate this policy using 500 out-of-sample scenarios. In the SDDP.jl package, this procedure is called policy simulation. Specifically, we chose the `SDDP.InSampleMonteCarlo' simulation scheme. SDDP.jl package considers this simulation scheme as in-sample because it considers the same random variable distribution as the training procedure.

The two-stage approximation we considered is a restriction of the multi-stage SP model. Finally, we create out-of-sample validation scenarios to test the performance of decision policies obtained from both models.

\section{Numerical Experiments and Sensitivity Analysis}\label{sec:Results}

The goal of our numerical experiments is twofold: ($i$) to demonstrate the value of multi-stage SP with respect to their two-stage counterpart and the MV problem; ($ii$) to evaluate the impact of different modeling and operational choices on the system performance. Specifically, we conduct sensitivity analysis regarding the granularity of decision stages, initial inventory level, length of the bale sequencing pattern, and the ordering in the bale sequence. 

We implemented every optimization model considered in the paper in Julia 1.6.1 using the mathematical optimization modeling package JuMP~\citep{DunningHuchetteLubin2017}. We conducted the experiments on Clemson University's high-performance computing cluster, the Palmetto Cluster, and used 16 nodes with 32 GB RAM. We solve the multi-stage SP model using the SDDP.jl package~\citep{dowson_sddp.jl}. SDDP.jl is a Julia package for solving multi-stage convex stochastic programming problems using SDDP. We use the commercial solver Gurobi as the optimization solver within the package.

\subsection{Data Description}\label{sec:Case}
We develop a case study using historical data about the characteristics of biomass and the performance of equipment of PDU~\citep{Gulcan2021}. We also consulted the experts and operators of the PDU during model development, verification, and validation. Our data set summarizes the sensors reading of moisture content collected during the processing of switchgrass bales at PDU. The moisture content of bales in our case study varies from $3\%$ to $30\%$, and bales are grouped into low ($3\%$ to $12\%$), medium ($12\%$ to $20\%$), and high ($20\%$ to $30\%$) moisture levels. We assume that the probability distribution of moisture content within each level follows a uniform distribution. 

Work by~\citep{HANSEN2019BiomassSC} provides additional data about the bulk density of switchgrass when different harvesting equipment are used. The average of the dry biomass densities that they report is $203.04\ kg/m^3$. Note that even though we use the average dry density value, the wet density and the biomass flow vary due to the random moisture content. 

Finally, we use data generated by the Discrete Element Method (DEM) model to create regression functions \eqref{reg1} and \eqref{reg2}. These equations represent the relationship between moisture content, particle size distribution, and bulk density.  
Regressions \eqref{reg1} and \eqref{reg2} present biomass density after processed at grinder 1 and grinder 2, respectively.

\begin{equation}\label{reg1}
\tilde{d}_{g_1,t} = 56.183 + 65.312 \tilde{m}_{g_1,t}  - 8.473\rho^{50}_{g_1} \\
\end{equation}
\begin{equation}\label{reg2}
\tilde{d}_{g_2,t} = 186.348 + 206.1697 \tilde{m}_{g_2,t} - 110.302 \rho^{50}_{g_2} \\
\end{equation}
where $d_{it}$ represents the biomass density (dependent variable), $m_{it}$ represents the moisture level of the biomass processed in equipment $i = g_1, g_2$ (indices $g_1$ and $g_2$ correspond to grinder 1 and grinder 2, respectively) and $\rho^{50}_i$ represents the $50$-th percentile of the particle size distribution of biomass processed in equipment $i$. The coefficients $\{(\alpha^0_i,\alpha^1_i,\alpha^2_i)\}_{i=g_1,g_2}$ in the regression function are determined based on the DEM by \citep{guo2020discrete} that simulates the process in grinders 1 and 2. Tables~\ref{tab:PSD} to \ref{tab:processing_changes} in Appendix~\ref{App:DataTables} summarize the input data used. 

\subsection{Selection of Sample Sizes for SP Models}\label{sec:ScenarioNum}
We first conducted in-sample and out-of-sample stability tests to identify the appropriate number of scenarios to use in the policy training and policy evaluation phases for the SP models under consideration. Let $N_t$ be the number of realizations used in each stage for solving the multi-stage SP model via the SDDP algorithm, and let $S_S$ be the number of sample paths used for solving the two-stage SP model. Additionally, we let $S_V$ be the number of sample paths created to test the performance of any decision policy (provided by either the multi-stage SP, two-stage SP or MV problem) in the out-of-sample evaluation phase. 

\paragraph{Multi-stage SP Model - $N_t$}
Table~\ref{tbl:N_t} summarizes the optimal objective value information from the multi-stage SP model with the number of realizations per stage ($N_t$) in the training sample ranging from $250$ to $750$. Recall that the SDDP algorithm gives a deterministic lower bound (\textbf{LB}) and a statistical upper bound (\textbf{CI}) associated with the obtained decision policy is given by an out-of-sample experiment. We used $500$ out-of-sample scenarios in these experiments. The results of Table~\ref{tbl:N_t} show that there is no significant difference in the out-of-sample performance of the policy for different values of $N_t$ between $250$ and $750$. Thus, we chose the middle value and set $N_t\ =\ 500$ for our experiments. 

\paragraph{Two-stage SP Model - $S_S$}
Table~\ref{tbl:S_S} presents the optimal objective value of the two-stage approximation with the number of (training) scenarios $S_S$ ranging from $500$ to $1,500$. Table~\ref{tbl:S_S} presents the in-sample performance of the two-stage model.  Again, we did not observe any significant difference in the optimal objective value for the values of $S_S$ within the range between $500$ and $1,500$. Thus, we set $S_S\ =\ 1,000$ for our experiments.

\begin{table}[htbp]
\centering
	\scriptsize
	{%
	\caption{In-sample stability test - \# of scenarios used (per stage) in the multi-stage SP model.}
	\label{tbl:N_t}
		 \setlength{\tabcolsep}{0.2em}
        	\begin{tabular}{lccc}
			\toprule
			& \multicolumn{3}{c}{\textbf{Number of Realizations}
			$\boldsymbol{N_t}$}\\
			\cmidrule(lr){2-4}
			&	250 	&	500* 	&	750\\	
			\cmidrule(lr){2-2}
			\cmidrule(lr){3-3}
			\cmidrule(lr){4-4}
			\textbf{CI* (\$)}	 & $\left \lbrack	124.26	- 124.55 \right\rbrack$  & $\left \lbrack	124.41 - 	124.69 \right\rbrack$ & $\left \lbrack	124.38 -	124.66\right\rbrack$ \\
            \textbf{LB (\$)}	&	124.36 &	124.48	&	124.45	\\
			\bottomrule\\
			\multicolumn{4}{l}{*$95\%$ Confidence Interval}\\
	\end{tabular}}
\end{table}%

\begin{table}[htbp]
\centering
	\scriptsize
	{%
	\caption{In-sample stability test - \# of scenarios used in the two-stage SP model.}
	\label{tbl:S_S}
        \begin{tabular}{lccc}
			\toprule
			& \multicolumn{3}{c}{\textbf{Number of Scenarios} $\boldsymbol{S_S}$}\\
			\cmidrule(lr){2-4}
             & 500 & 1,000* & 1,500\\
            \cmidrule(lr){2-2}
			\cmidrule(lr){3-3}
			\cmidrule(lr){4-4}
            \textbf{Objective (\$)} &	128.25 &	128.28 & 128.32\\
			\bottomrule
			\\
	\end{tabular}}
\end{table}%

\paragraph{Out-of-sample Evaluation - $S_V$}
Table~\ref{tbl:S_V} presents the CI on the optimal cost for different models: multi-stage SP, two-stage SP, and a deterministic model, under different number of ($S_V$) out-of-sample scenarios (ranging from $250$ to $750$). For these experiments, we use decision policies trained by multi-stage SP and two-stage SP models with $N_t\ = 500$ and $S_S\ =\ 1,000$, respectively. Based on the out-of-sample stability tests, we set $S_V\ =\ 500$ for our experiments.

\begin{table}[htbp]
\centering
	\scriptsize
	{%
	\caption{Out-of-sample stability test - \# of sample paths used for performance evaluation.}
	\label{tbl:S_V}
        \begin{tabular}{lcc}
        			\toprule
        		   \textbf{Model} & $\boldsymbol{S_V}$ & \textbf{CI* (\$)} \\
        		   \cmidrule(lr){1-1}
        		    \cmidrule(lr){2-2}
        			\cmidrule(lr){3-3}
        			Multi-stage & \multirow{3}{*}{250} & $\left \lbrack 124.28 -	124.68	\right\rbrack$ \\
        			MV Problem & & $\left \lbrack 130.09 -	130.59	\right\rbrack$\\
        			Two-stage& & $\left \lbrack 128.10 -	128.47	\right\rbrack$ \\
        			\midrule
        			Multi-stage & \multirow{3}{*}{500*} & $\left \lbrack 124.41 -	124.69	\right\rbrack$ \\
        			MV Problem & &  $\left \lbrack 130.21 -	130.56	\right\rbrack$\\
        			Two-stage& & $\left \lbrack 128.22 -	128.48 \right\rbrack$ \\
        			\midrule
        			Multi-stage & \multirow{3}{*}{750} & $\left \lbrack 124.41 -	124.64	\right\rbrack$ \\
        			MV Problem & &  $\left \lbrack 130.22 -	130.53	\right\rbrack$\\
        			Two-stage& & $\left \lbrack 128.23 -	128.44 \right\rbrack$ \\
        			\bottomrule\\
        			\multicolumn{3}{l}{*$95\%$ Confidence Interval}
        	\end{tabular}}
\end{table}%

\subsection{Performance of the Multi-stage SP Model}\label{sec:model-evaluation}
In this section, we demonstrate the value of the proposed multi-stage SP model by comparing it with the two-stage approximation and the deterministic MV problem that we presented in Sections~\ref{sec:2stage} and \ref{sec:Deterministic}, respectively.

\subsubsection{Base-case problem} \label{sec:base-analyis}
 First, we define a base-case problem, which considers a planning horizon of processing 50 bales, among which 60\% are low-moisture, 20\% are medium-moisture, and 20\% are high-moisture bales. We assume that the bales are sequenced based on a repeating pattern of one high moisture bale, one medium moisture bale, and three low moisture bales over the planning horizon. We call this type of sequencing pattern the `short' sequence. In the base-case problem, we set the target reactor feeding rate to be $2.95\ dt/hr$.

Table~\ref{tbl:Performance_base} summarizes the performance of the three models under consideration for the base-case problem. The results of Table~\ref{tbl:Performance_base} demonstrate the value of the proposed multi-stage SP model: it provides an operational decision policy that costs 3.1\% less than the static policy obtained by the two-stage SP model, and 4.7\% less than the static policy obtained by the MV problem, in terms of the total expected cost. Not surprisingly, this superiority of solution quality comes with a price of excessive computational time. Also, the CIs of these models do not intersect. Thus difference between solutions provided by each policy is statistically significant.

\begin{figure}[h]
    \centering
    \begin{subfigure}[t]{0.4\textwidth}
        \centering
        \includegraphics[width=\textwidth]{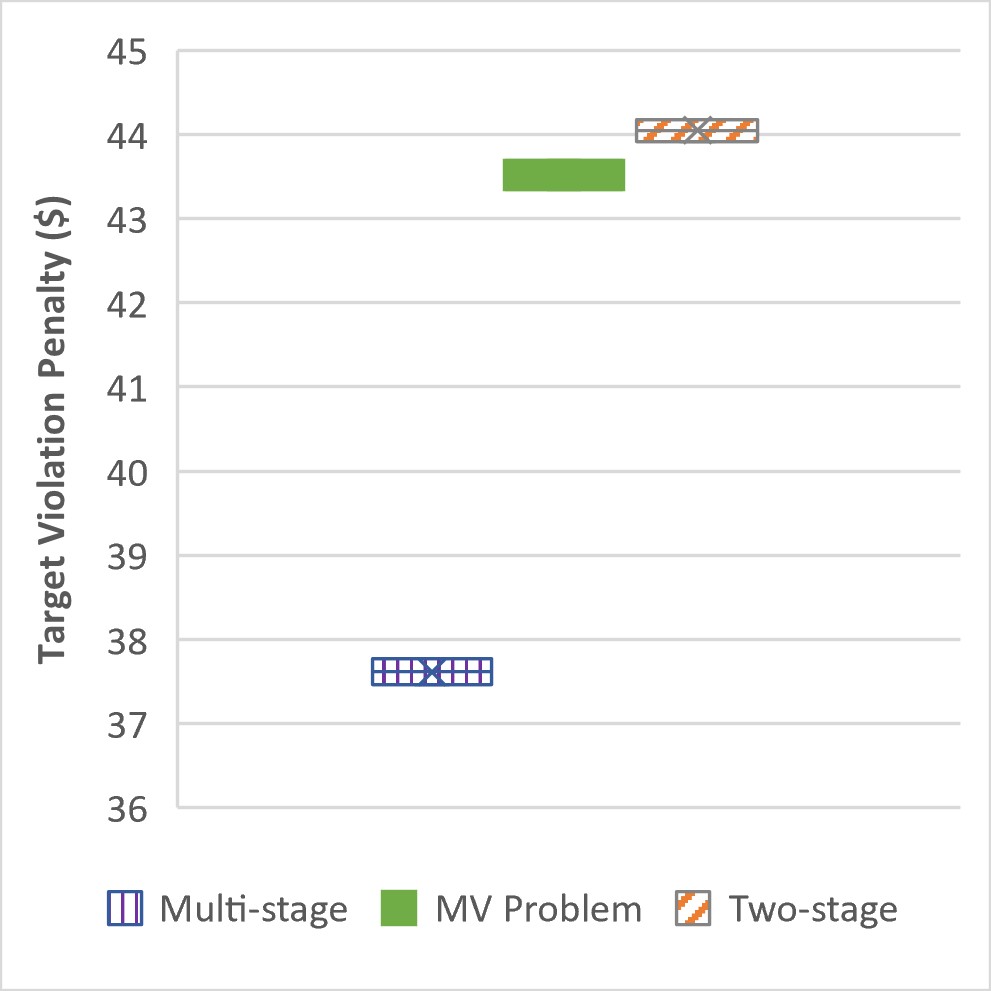}
    \caption{Target violation penalties for different models.}
    \label{fig:base-case}
    \end{subfigure}
    \hfill
    \begin{subfigure}[t]{0.55\textwidth}
        \centering
    	\includegraphics[width=\textwidth]{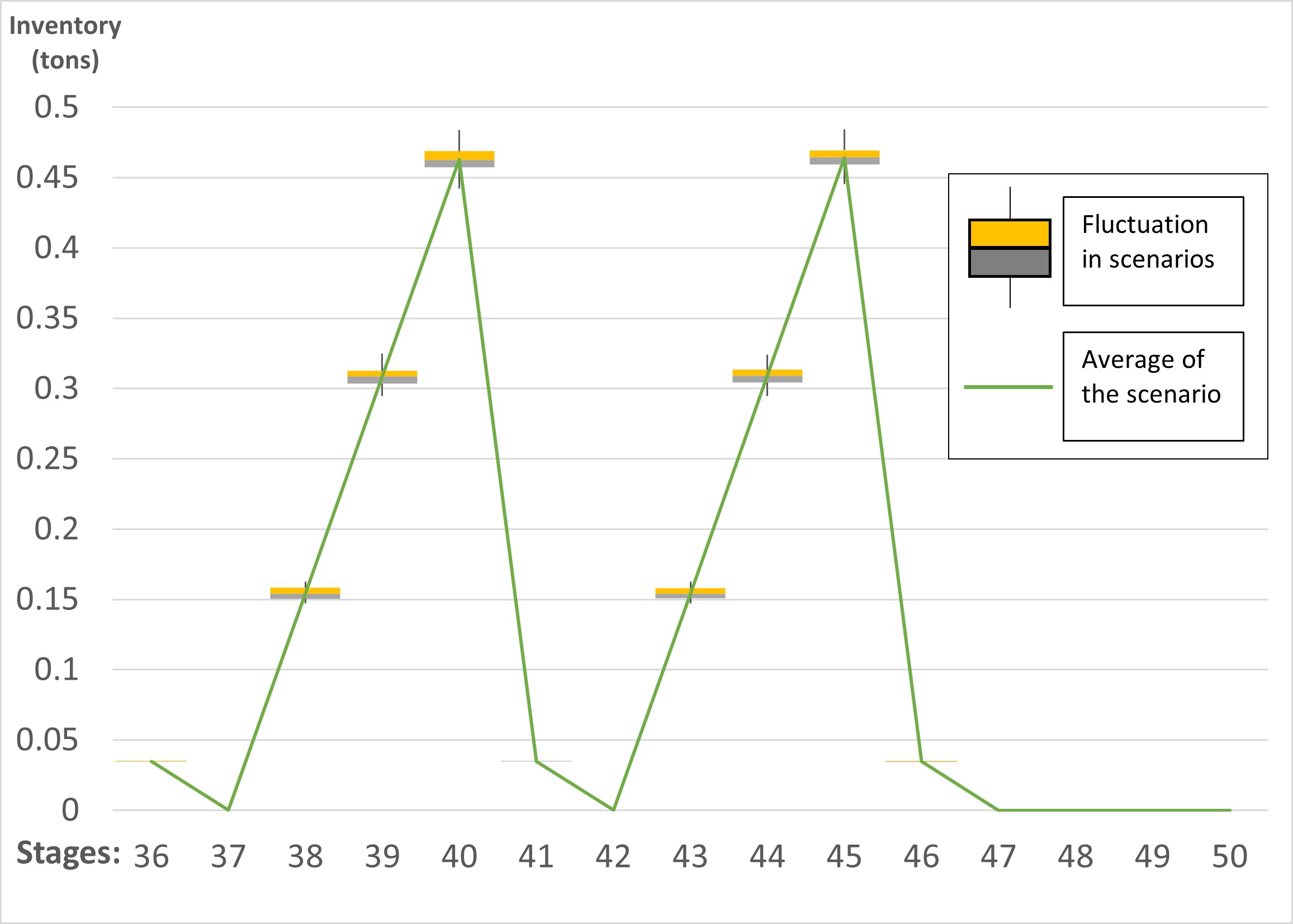}
    	\caption{Adaptive inventory levels given by the multi-stage SP.}
    	\label{fig:base-case_Multi_Inv}
    \end{subfigure}
    \caption{Model performance on the base-case problem.}
\end{figure}

As shown in Figure~\ref{fig:base-case}, the operational policy obtained using the multi-stage SP model performs better in terms of achieving the target reactor rate by adaptively setting the equipment processing speed based on the changes in the moisture content measured by the sensors. 

\begin{table}[htbp]
\centering
	\scriptsize
	{%
	\caption{Model performance on the base-case problem.}
        	\label{tbl:Performance_base}
        \begin{tabular}{lccc}
        			\toprule
        			& \multicolumn{3}{c}{\textbf{50 Bales}}\\
        			\cmidrule(lr){2-4}
        			& 
        			&  \textbf{Gap**} &  \textbf{Run Time} \\
        			\textbf{Model} 
        			& \textbf{CI* (\$)}
        			&  \textbf{(\%)} &  \textbf{(sec)} \\
        			\cmidrule(lr){1-1}
        			\cmidrule(lr){2-2}
        			\cmidrule(lr){3-3}
        			\cmidrule(lr){4-4}
        			Multi-stage
        			& $\left \lbrack 124.41 - 124.69 \right\rbrack$ &	- & 644\\
                    MV Problem 
                    & $\left \lbrack130.21 -	130.56\right\rbrack$	& 4.7 &	7\\
                    Two-stage 
                    & $\left \lbrack128.22 -	128.48\right\rbrack$ &	3.1	&74\\
        			\bottomrule\\
        			\multicolumn{4}{l}{*$95\%$ Confidence Interval}\\
        			\multicolumn{4}{l}{**Gap := Percentage difference compared to the multi-stage SP}\\
        			\multicolumn{4}{l}{using the mean performance (middle point of the CI).}
        	\end{tabular}}
\end{table}%

Figure~\ref{fig:base-case_Multi_Inv} shows the inventory levels given by the multi-stage SP model for stages 36 to 50. In Figure~\ref{fig:base-case_Multi_Inv}, the solid line represents the average inventory level over $500$ sample paths (for the out-of-sample evaluation) in each stage. We use box-and-whisker plots to represent the fluctuations in the inventory levels for different sample paths in each stage. As we see from Figure~\ref{fig:base-case_Multi_Inv}, the inventory level changes in a cyclic pattern every five stages where the moisture level in the bale sequence completes a whole pattern of one high, one medium, and three low. The inventory is accumulated while processing low-moisture bales because the system can run faster. The system uses the accumulated inventory while processing high- and medium-moisture bales. The decrease in the inventory level is steeper when high moisture bales are processed because the system runs slower and requires higher inventory to meet the target reactor rate. In the last three stages, we observe no inventory due to the end-of-horizon effect.

Table~\ref{tbl:Inv_comparison_base} compares inventory levels set by the two-stage SP and MV problem with the average inventory level from the multi-stage SP solutions. We see that the two-stage SP model results in lower inventory levels, which reduces the inventory holding costs, but in return, increases the penalty cost caused by the violation of the reactor's target rate, as shown in Figure~\ref{fig:base-case}.

\begin{table*}[htbp]
\centering
	\scriptsize
	{%
	\caption{Percentage difference in the inventory level obtained by the two-stage SP model and MV problem in stage 36 to stage 50 for the base-case problem compared to the multi-stage SP model.}
	\label{tbl:Inv_comparison_base}
		\setlength{\tabcolsep}{0.47em}
        \begin{tabular}{lccccccccccccccc}
        			\toprule
        			& \multicolumn{15}{c}{\textbf{Stages}}\\
        			\cmidrule(lr){2-16}
        			\textbf{Model} &	36 &	37 &	38 &	39	&40 &	41&	42&	43&	44&	45&	46&	47&	48&	49&	50\\
        			\cmidrule(lr){1-1}
        			\cmidrule(lr){2-2}
        			\cmidrule(lr){3-3}
        			\cmidrule(lr){4-4}
        			\cmidrule(lr){5-5}
        			\cmidrule(lr){6-6}
        			\cmidrule(lr){7-7}
        			\cmidrule(lr){8-8}
        			\cmidrule(lr){9-9}
        			\cmidrule(lr){10-10}
        			\cmidrule(lr){11-11}
        			\cmidrule(lr){12-12}
        			\cmidrule(lr){13-13}
        			\cmidrule(lr){14-14}
        			\cmidrule(lr){15-15}
        			\cmidrule(lr){16-16}
                    \textbf{Two-stage} &	-0.1\%& 	0.0\%& 	-3.6\%& -3.1\%& 	-2.7\%&	0.0\%&	0.0\%&	-3.7\%&	-3.3\%&	-3.1\%&	0.0\%&	0.0\%&	0.0\%&	0.0\%&	0.0\%\\
                    \textbf{MV Problem}&	0.4\%&	0.0\%&	0.1\%&	0.1\%&	0.1\%&	0.5\%&	0.0\%&	0.0\%&	0.0\%&	-0.2\%&	0.6\%&	0.0\%&	0.0\%&	0.0\%&	0.0\%\\
        			\bottomrule
        	\end{tabular}}
\end{table*}%

\subsubsection{Length of planning horizon}

The goal of the following experiments is to assess the value of the multi-stage SP as the number of stages increases.
In our problem setting, the length of the planning horizon is identified by the number of bales that the system will process during the operation time. 
Table~\ref{tbl:Performance_NumStage} summarizes the performance comparison of the three models on problems with different planning horizons. For a fair comparison, Table~\ref{tbl:Performance_NumStage} records the total cost per bale. In addition to the base-case problem of processing 50 bales, we considered problems with a planning horizon of $10,\ 25$, and $100$ stages/bales.

\begin{table}[htbp]
\centering
	\scriptsize
	{%
	\caption{Model performances for different length of planning horizons.}
        	\label{tbl:Performance_NumStage}
        \begin{tabular}{cccc}
        			\toprule
        			& \multicolumn{3}{c}{\textbf{10 Bales}}\\
        			\cmidrule(lr){2-4}
        			 & \textbf{Cost per Bale}
        			&  \textbf{Gap*} &  \textbf{Run Time} \\
        			\textbf{Model} & \textbf{CI (\$)}
        			&  \textbf{(\%)} &  \textbf{(sec)} \\
        			\cmidrule(lr){1-1}
        			\cmidrule(lr){2-2}
        			\cmidrule(lr){3-3}
        			\cmidrule(lr){4-4}
        			Multi-stage
        			& $\left \lbrack 3.389 - 3.399 \right\rbrack$ &	- & 128 \\
        			MV Problem 
                    & $\left \lbrack3.453 -	3.465\right\rbrack$ &	1.9 & 7\\
                    Two-stage 
                    & $\left \lbrack3.431 -	3.440\right\rbrack$ &	1.2	&18	\\
                    \midrule
                    & \multicolumn{3}{c}{\textbf{25 Bales}}\\
                    \cmidrule(lr){2-4}
        			 & \textbf{Cost per Bale}
        			&  \textbf{Gap*} &  \textbf{Run Time} \\
        			\textbf{Model} & \textbf{CI (\$)}
        			&  \textbf{(\%)} &  \textbf{(sec)} \\
        			\cmidrule(lr){1-1}
        			\cmidrule(lr){2-2}
        			\cmidrule(lr){3-3}
        			\cmidrule(lr){4-4}
        			Multi-stage & $\left \lbrack 2.711 - 2.719 \right\rbrack$ &	- & 300	\\
                    MV Problem 
                    & $\left \lbrack2.816 -	2.826\right\rbrack$	&3.9	&7\\
        			Two-stage 
                    & $\left \lbrack2.779 -	2.786\right\rbrack$ &	2.5	&38\\
                    \midrule
        			& \multicolumn{3}{c}{\textbf{100 Bales}}\\
        			\cmidrule(lr){2-4}
        			 & \textbf{Cost per Bale}
        			&  \textbf{Gap**} &  \textbf{Run Time} \\
        			\textbf{Model} & \textbf{CI* (\$)}
        			&  \textbf{(\%)} &  \textbf{(sec)} \\
        			\cmidrule(lr){1-1}
        			\cmidrule(lr){2-2}
        			\cmidrule(lr){3-3}
        			\cmidrule(lr){4-4}
        			Multi-stage & $\left \lbrack 2.375 - 2.380	\right\rbrack$ &	- &	2641\\
                    MV Problem
                    & $\left \lbrack2.498 -	2.503\right\rbrack$	&5.2 &11\\
                    Two-stage
                    & $\left \lbrack2.456 -	2.460\right\rbrack$ &	3.4	&299\\
        			\bottomrule\\
        			\multicolumn{4}{l}{*$95\%$ Confidence Interval}\\
        			\multicolumn{4}{l}{**Gap := Percentage difference compared to the multi-stage SP}\\
        			\multicolumn{4}{l}{using the mean performance (middle point of the CI).}
        	\end{tabular}}
\end{table}%

The results of Tables \ref{tbl:Performance_base} and \ref{tbl:Performance_NumStage} indicate that as the number of stages increases, the value of using a multi-stage stochastic program increases. Our model proposes that bales are processed based on a sequence that repeats itself. Thus, one might consider using a `truncate-and-repeat' strategy by obtaining an optimal policy from solving a problem with a small number of stages and then repeatedly using the resulting policy as a heuristic approach over the entire planning horizon. Table~\ref{tbl:RepeatingStrategy} summarizes the outcome of applying this strategy. We can see that for our problem, this `truncate-and-repeat' strategy results in a worse performance that is statistically significant compared to the policy obtained from solving the full-length problem. This is attributed to the `end-of-horizon' effect in multi-stage SP problems~\citep{Shapiro_2011}. If the time is limited and the decision needs to be made quickly, rather than using this `truncate-and-repeat' strategy one can use the MV problem or two-stage model. Tables \ref{tbl:Performance_base} and \ref{tbl:Performance_NumStage} show that both models provide higher-quality solutions than the solutions obtained by the `truncate-and-repeat' strategy and require much shorter computational time. 

\begin{table}[htbp]
\centering
	\scriptsize
	{%
	\caption{The performance of the `truncate-and-repeat' strategy.}
        	\label{tbl:RepeatingStrategy}
		\setlength{\tabcolsep}{0.5em}
        \begin{tabular}{cccc}
        			\toprule
        			& \textbf{5 Times} & \textbf{2 Times} & \textbf{1 Time}\\
        			& \textbf{10-Bale Policy} & \textbf{25-Bale Policy} & \textbf{50-Bale Policy}\\
        		    \cmidrule(lr){2-2}
        			\cmidrule(lr){3-3}
        			\cmidrule(lr){4-4}
                   \textbf{Total Cost (\$)} & $\left \lbrack 169.5 - 169.9
                	\right\rbrack$ & $\left \lbrack 135.6 -	135.9
                	\right\rbrack$ & $\left \lbrack 124.4 -	124.7	\right\rbrack$\\
                	\textbf{Cost per }& \multirow{2}{*}{$\left \lbrack 3.389 -	3.399
                	\right\rbrack$}
                	& \multirow{2}{*}{$\left \lbrack 2.711 - 2.719
                	\right\rbrack$}
                	&\multirow{2}{*}{$\left \lbrack 2.488 - 2.494
	                \right\rbrack$}\\
	                \textbf{Bale (\$)} \\
                    \textbf{Gap* (\%)} & 36 & 9 & - \\
        			\bottomrule\\
        			\multicolumn{4}{l}{*Gap := Percentage difference compared to base-case of 50 bales}\\
        			\multicolumn{4}{l}{using the mean performance (middle point of the CI).}
        	\end{tabular}}
\end{table}%

\subsubsection{The mixing ratio of biomass bales based on their moisture levels}

The mix of biomass bales based on their moisture levels can greatly impact the performance of biomass processing system~\citep{Gulcan2021}. The goal of the following experiments is to analyze how the value of multi-stage SP varies under different mixes of biomass bales. In addition to the base-case problem of processing 60\% low - 20\% med - 20\% high moisture bales, we considered processing all low, all high, 60\% med - 20\% low - 20\% high and 60\% high - 20\% med - 20\% low moisture bales.  For the bale mixes 60\% med - 20\% low - 20\% high and 60\% high - 20\% med - 20\% low, we considered the bales are sequenced using the short sequencing strategy with internal patterns of one high moisture bale, three medium moisture bales, and one low moisture bale, and three high moisture bales, one medium moisture bale, and one low moisture bale, respectively. For each mix, we experimented with three different target reactor rates; $2.50\ dt/hr$, $2.72\ dt/hr$, and $2.95\ dt/hr$. We considered that the remainder of the problem parameters is the same as the base-case problem.

\begin{table}[htbp]
\centering
	\scriptsize
	{%
	\caption{Model performances for different mixing ratios of bales.}}
		
		\begin{subtable}[t]{0.38\textwidth}
		\setlength{\tabcolsep}{0.2em}
		\centering
		\scriptsize	{%
		\caption{Target rate $2.95\ dt/hr$.}
	    \label{tbl:MixingRatio2.95}
		\begin{tabular}{lcc}
			\toprule
			\multicolumn{3}{c}{\textbf{All Low}}\\
			\cmidrule(lr){1-3}
			& 	&  \textbf{Gap**}  \\
			\textbf{Model} & \textbf{CI* (\$)} &  \textbf{(\%)}  \\
			\cmidrule(lr){1-1}
			\cmidrule(lr){2-2}
			\cmidrule(lr){3-3}
			Multi-stage
			& $\left \lbrack 0.00 - 0.00\right\rbrack$ &	-  \\
			MV Problem 
			& $\left \lbrack0.00 - 0.00\right\rbrack$ &	0.0 \\
			Two-stage 
			& $\left \lbrack0.00 - 0.00\right\rbrack$ &	0.0	\\
			\midrule
			\multicolumn{3}{c}{ $\boldsymbol{60\%\ Low$-$20\%\ Med$-$20\%\ High}$}\\
			\cmidrule(lr){1-3}
			& 	&  \textbf{Gap**}  \\
			\textbf{Model} & \textbf{CI* (\$)} &  \textbf{(\%)} \\
			\cmidrule(lr){1-1}
			\cmidrule(lr){2-2}
			\cmidrule(lr){3-3}
			Multi-stage
			& $\left \lbrack 124.41 - 124.69 \right\rbrack$ &	- \\
            MV Problem 
            & $\left \lbrack130.21 -	130.56\right\rbrack$	& 4.7\\
            Two-stage 
            & $\left \lbrack128.22 -	128.48\right\rbrack$ &	3.1\\
			\midrule
			\multicolumn{3}{c}{ $\boldsymbol{60\%\ Med$-$20\%\ Low$-$20\%\ High}$}\\
			\cmidrule(lr){1-3}
			& 	&  \textbf{Gap**}  \\
			\textbf{Model} & \textbf{CI* (\$)} &  \textbf{(\%)} \\
			\cmidrule(lr){1-1}
			\cmidrule(lr){2-2}
			\cmidrule(lr){3-3}
			Multi-stage & $\left \lbrack 225.91 - 226.00	\right\rbrack$ &	- \\
			MV Problem
			& $\left \lbrack228.14 - 228.29\right\rbrack$	&1.0 \\
			Two-stage
			& $\left \lbrack227.78 - 227.86\right\rbrack$ &	0.8\\
			\midrule
			\multicolumn{3}{c}{ $\boldsymbol{60\%\ High$-$20\%\ Med$-$20\%\ Low}$}\\
			\cmidrule(lr){1-3}
			& 	&  \textbf{Gap**}  \\
			\textbf{Model} & \textbf{CI* (\$)} &  \textbf{(\%)} \\
			\cmidrule(lr){1-1}
			\cmidrule(lr){2-2}
			\cmidrule(lr){3-3}
			Multi-stage & $\left \lbrack 591.29 - 591.45	\right\rbrack$ &	- \\
			MV Problem
			& $\left \lbrack593.24 - 593.44\right\rbrack$	&0.3\\
			Two-stage
			& $\left \lbrack592.97 - 593.13\right\rbrack$ &	0.3\\
			\midrule
			\multicolumn{3}{c}{\textbf{All High}}\\
			\cmidrule(lr){1-3}
			& 	&  \textbf{Gap**}  \\
			\textbf{Model} & \textbf{CI* (\$)} &  \textbf{(\%)} \\
			\cmidrule(lr){1-1}
			\cmidrule(lr){2-2}
			\cmidrule(lr){3-3}
			Multi-stage
			& $\left \lbrack 1033.28	 - 1033.47\right\rbrack$ &	- \\
			MV Problem 
			& $\left \lbrack1033.28	 - 1033.47\right\rbrack$ &	0.0\\
			Two-stage 
			& $\left \lbrack1033.28	 - 1033.47\right\rbrack$ &	0.0	\\
			\bottomrule\\
			\multicolumn{3}{l}{*$95\%$ Confidence Interval}\\
			\multicolumn{3}{l}{**Gap := Average change in objective}\\
			\multicolumn{3}{l}{values compared to multi-stage SP model.}
	\end{tabular}
	}
	\end{subtable}
	\hfill
	\begin{subtable}[t]{0.3\textwidth}
	\centering
	\scriptsize{%
	\caption{Target rate $2.72\ dt/hr$.}
	        \label{tbl:MixingRatio2.72}
		\begin{tabular}{cc}
			\toprule
			\multicolumn{2}{c}{\textbf{All Low}}\\
			\cmidrule(lr){1-2}
			&  \textbf{Gap**}  \\
		    \textbf{CI* (\$)} &  \textbf{(\%)} \\
			\cmidrule(lr){1-1}
			\cmidrule(lr){2-2}
			$\left \lbrack 0.00 - 0.00\right\rbrack$ &	- \\
			$\left \lbrack0.00 - 0.00\right\rbrack$ &	0.0\\
			$\left \lbrack0.00 - 0.00\right\rbrack$ &	0.0\\
			\midrule
			\multicolumn{2}{c}{$\boldsymbol{60\%\ Low$-$20\%\ Med$-$20\%\ High}$}\\
			\cmidrule(lr){1-2}
			&  \textbf{Gap**}  \\
			\textbf{CI* (\$)} &  \textbf{(\%)} \\
			\cmidrule(lr){1-1}
			\cmidrule(lr){2-2}
			$\left \lbrack 68.08 -	68.20 \right\rbrack$ &	- \\
            $\left \lbrack70.27 -	70.49\right\rbrack$	& 3.3\\
            $\left \lbrack68.84 -	69.00\right\rbrack$ &	1.1\\
			\midrule
			\multicolumn{2}{c}{$\boldsymbol{60\%\ Med$-$20\%\ Low$-$20\%\ High}$}\\
			\cmidrule(lr){1-2}
			&  \textbf{Gap**}  \\
			 \textbf{CI* (\$)} &  \textbf{(\%)} \\
			\cmidrule(lr){1-1}
			\cmidrule(lr){2-2}
			$\left \lbrack 97.70 - 97.8\right\rbrack$ &	-\\
			$\left \lbrack101.48 - 101.70\right\rbrack$	&3.9\\
			$\left \lbrack100.57 - 100.72\right\rbrack$ &	2.9\\
			\midrule
			\multicolumn{2}{c}{$\boldsymbol{60\%\ High$-$20\%\ Med$-$20\ \%Low}$}\\
			\cmidrule(lr){1-2}
			&  \textbf{Gap**}  \\
			\textbf{CI* (\$)} &  \textbf{(\%)} \\
			\cmidrule(lr){1-1}
			\cmidrule(lr){2-2}
			$\left \lbrack 431.27	- 431.46	\right\rbrack$ &	-\\
			$\left \lbrack434.16 - 434.39\right\rbrack$	&0.7\\
			$\left \lbrack433.76 - 433.93\right\rbrack$ &	0.6\\
			\midrule
			\multicolumn{2}{c}{\textbf{All High}}\\
			\cmidrule(lr){1-2}
			& 	\textbf{Gap**}  \\
			\textbf{CI* (\$)} &  \textbf{(\%)} \\
			\cmidrule(lr){1-1}
			\cmidrule(lr){2-2}
			$\left \lbrack 825.04 - 825.22\right\rbrack$ &	- \\
			$\left \lbrack825.04 -	825.22\right\rbrack$ &	0.0 \\
			$\left \lbrack825.04 -	825.22\right\rbrack$ &	0.0	\\
			\bottomrule\\
			\multicolumn{2}{l}{*$95\%$ Confidence Interval}\\
			\multicolumn{2}{l}{**Gap := Average change in}\\
			\multicolumn{2}{l}{objective values compared}\\
			\multicolumn{2}{l}{to multi-stage SP model.}
	\end{tabular}
	}
	\end{subtable}
	\hfill
	\begin{subtable}[t]{0.3\textwidth}
	\centering
	\scriptsize{%
	\caption{Target rate $2.50\ dt/hr$.}
	        \label{tbl:MixingRatio2.5}
		\begin{tabular}{cc}
			\toprule
			\multicolumn{2}{c}{\textbf{All Low}}\\
			\cmidrule(lr){1-2}
			&  \textbf{Gap**}  \\
		    \textbf{CI* (\$)} &  \textbf{(\%)} \\
			\cmidrule(lr){1-1}
			\cmidrule(lr){2-2}
			$\left \lbrack 0.00 - 0.00\right\rbrack$ &	- \\
			$\left \lbrack0.00 - 0.00\right\rbrack$ &	0.0\\
			$\left \lbrack0.00 - 0.00\right\rbrack$ &	0.0\\
			\midrule
			\multicolumn{2}{c}{$\boldsymbol{60\%Low$-$20\%Med$-$20\%High}$}\\
			\cmidrule(lr){1-2}
			&  \textbf{Gap**}  \\
			\textbf{CI* (\$)} &  \textbf{(\%)} \\
			\cmidrule(lr){1-1}
			\cmidrule(lr){2-2}
			$\left \lbrack 68.08 -	68.20 \right\rbrack$ &	-\\
            $\left \lbrack70.27 -	70.49\right\rbrack$	& 3.3\\
            $\left \lbrack68.84 -	69.00\right\rbrack$ &	1.1\\
			\midrule
			\multicolumn{2}{c}{$\boldsymbol{60\%Med$-$20\%Low$-$20\%High}$}\\
			\cmidrule(lr){1-2}
			&  \textbf{Gap**}  \\
			 \textbf{CI* (\$)} &  \textbf{(\%)} \\
			\cmidrule(lr){1-1}
			\cmidrule(lr){2-2}
			$\left \lbrack 43.22 -	43.31\right\rbrack$ &	-\\
			$\left \lbrack45.64 -	45.86\right\rbrack$	&5.7\\
			$\left \lbrack43.79 -	43.86\right\rbrack$ &	1.3\\
			\midrule
			\multicolumn{2}{c}{$\boldsymbol{60\%High$-$20\%Med$-$20\%Low}$}\\
			\cmidrule(lr){1-2}
			&  \textbf{Gap**}  \\
			\textbf{CI* (\$)} &  \textbf{(\%)} \\
			\cmidrule(lr){1-1}
			\cmidrule(lr){2-2}
			$\left \lbrack 284.99 -	285.16	\right\rbrack$ &	-\\
			$\left \lbrack289.65 -	289.90\right\rbrack$	&1.6\\
			$\left \lbrack288.16 -	288.34\right\rbrack$ &	1.1\\
			\midrule
			\multicolumn{2}{c}{\textbf{All High}}\\
			\cmidrule(lr){1-2}
			& 	\textbf{Gap**}  \\
			\textbf{CI* (\$)} &  \textbf{(\%)} \\
			\cmidrule(lr){1-1}
			\cmidrule(lr){2-2}
			$\left \lbrack 616.77	- 616.95 \right\rbrack$ &	- \\
			$\left \lbrack616.77 -	616.95\right\rbrack$ &	0.0\\
			$\left \lbrack616.77 -	616.95\right\rbrack$ &	0.0\\
			\bottomrule\\
			\multicolumn{2}{l}{*$95\%$ Confidence Interval}\\
			\multicolumn{2}{l}{**Gap := Average change in}\\
			\multicolumn{2}{l}{objective values compared}\\
			\multicolumn{2}{l}{to multi-stage SP model.}
	\end{tabular}
	}
	\end{subtable}
\end{table}%

Tables~\ref{tbl:MixingRatio2.95}, \ref{tbl:MixingRatio2.72}, and \ref{tbl:MixingRatio2.5} summarize the results of processing different mixing ratios under different target reactor rates.  Tables~\ref{tbl:MixingRatio2.95}, \ref{tbl:MixingRatio2.72}, and \ref{tbl:MixingRatio2.5} show that operational costs increase as the ratio of higher moisture level bales increases in the mix. The reason is that when the system processes high moisture bales, equipment processing rates are lowered to prevent equipment clogging. 

The implications of Tables~\ref{tbl:MixingRatio2.95}, \ref{tbl:MixingRatio2.72}, and \ref{tbl:MixingRatio2.5} on the value of multi-stage SP are two-fold. First, when every bale has the same moisture level, there is no value to use a stochastic program to determine the inventory level. This is because the trivial solution of zero inventory is the optimal solution for these cases. When every bale has low moisture level, the target can be achieved without requiring an initial inventory and variation of realizations of moisture content only requires updates to the local control variables (equipment processing rates). Again in the case when moisture level is high, the best strategy is to maintain high feeding rates of the reactor. There are no expectations to accumulate inventory  when these bales are being processed. The infeed is the slowest under high moisture bales. Thus, if there is not enough inventory, the target reactor feeding rate cannot be reached. The costs presented in Tables~\ref{tbl:MixingRatio2.95}, \ref{tbl:MixingRatio2.72}, and \ref{tbl:MixingRatio2.5} incur solely due to the penalty for not meeting the targeted feeding rate. 

Second, when the system processes a mixture of bales, then, using a multi-stage SP model to optimize the biorefinery operations becomes valuable. The confidence intervals for these different models do not overlap, showing that the differences in different models' performance are statistically significant. Results of these tables show that different bale mix distributions and different target rates can lead to different values of multi-stage SP.

\subsection{Sensitivity Analysis: Granularity of Decision Stages} \label{sec:Sensitivity_Granularity}
In this section, we analyze the trade-offs among the granularity of decision stages, policy performance, and computational efforts. One of the modeling decisions that impacts the level of adaptability is the granularity of decision stages in the multi-stage SP. Recall that we have assumed thus far that the system processes a single bale at each stage. We now relax this assumption and consider two alternative stage definitions: ($i$) `Combined' stage definition, where a stage corresponds to the period in which the system will consecutively process biomass bales of the same moisture level; and ($ii$) `Detailed' stage definition, where a stage is the period in which the system processes one-third of a biomass bale. We assume that these new stage definitions do not conflict with the stage-wise independence assumption, and apply them on the base-case problem presented in Section~\ref{sec:model-evaluation}.

For these experiments, we update the objective function~\eqref{eqn:multi-stage-obj} of stage $t$ problem as follows:

\[
\min \  \beta_t(\mathbf{c}^{h \top} \mathbf{I}_t + c^pp_t)  + \mathcal{Q}_{t+1}(\mathbf{I}_{t}),
\]
where $\beta_t > 0$ is the number or fraction of bales that the system will process in stage $t$. The presence of $\beta_t$ ensures a unified total cost calculation under different stage definitions.

The `Combined' stage definition is motivated by the argument that we can use a similar operational decision policy for processing biomass bales of the same moisture level. On the other hand, the `Detailed' stage definition aims to interact with the system more frequently to increase adaptability to varying moisture content within the same bale. However, this increases the number of stages and thus increases the solution time for the multi-stage SP model.  Table~\ref{tbl:Granularity} summarizes the results of these two alternative stage definitions and compares them with the original stage definition.

\begin{table}[htbp]
	\centering
	\scriptsize
	{%
	\caption{Sensitivity analysis for granularity of decision stages.}
	\label{tbl:Granularity}
		\setlength{\tabcolsep}{0.5em}
		\begin{tabular}{lcccc}
			\toprule
			& 
			&  \textbf{Obj} & &  \textbf{Run Time} \\
			\textbf{Stage} & 
			&  \textbf{Change**} &  \textbf{Run Time}  &  \textbf{Change**} \\
			\textbf{Definition} 
			& \textbf{CI* (\$)}
			&  \textbf{(\%)} &  \textbf{(sec)} &  \textbf{(\%)}  \\
			\cmidrule(lr){1-1}
			\cmidrule(lr){2-2}
			\cmidrule(lr){3-3}
			\cmidrule(lr){4-4}
			\cmidrule(lr){5-5}
			Combined
			& $\left \lbrack166.09	- 166.38 \right\rbrack$	& 33.5 &	461 & -28.4\\
			Detailed 
			& $\left \lbrack123.21 -	123.36 \right\rbrack$ &	-1.0	& 2046& 217.7\\
			\bottomrule\\
			\multicolumn{5}{l}{*$95\%$ Confidence Interval}\\
			\multicolumn{5}{l}{**Change := Percentage change compared to the base-case with} \\
			\multicolumn{5}{l}{multi-stage SP in Table~\ref{tbl:Performance_base} using the mean performance}\\
			\multicolumn{5}{l}{(middle point of CI). }
	\end{tabular}}
\end{table}%

The results of Table~\ref{tbl:Granularity} show that the value of increasing the adaptability (with the detailed stage definition) is marginal in terms of solution quality (1\% improvement is given by the `Detailed' stage definition). However, this marginal gain is achieved at the expense of a 217.7\% increase in computational time. On the other hand, the `Combined' stage definition resulted in a 28.4\% reduction in computational time at the cost of a 33.5\% deterioration in the solution quality. 

\subsection{Sensitivity Analysis: Initial Inventory Level} \label{sec:Sensitivity_InitialInv}
In this section, we analyze the impact of the initial inventory level on costs and reactor utilization (i.e., target rate). Inventory holding is a common strategy to ensure high system reliability. We consider three different initial inventory levels for the storage equipment: empty, half-full, and full capacity, according to the common practice implemented at the PDU. For each of these initial inventory levels, we conducted experiments with different target reactor rates. 

\begin{figure}[h]
    \centering
    \begin{subfigure}[t]{0.48\textwidth}
        \centering
    	\includegraphics[width=\textwidth]{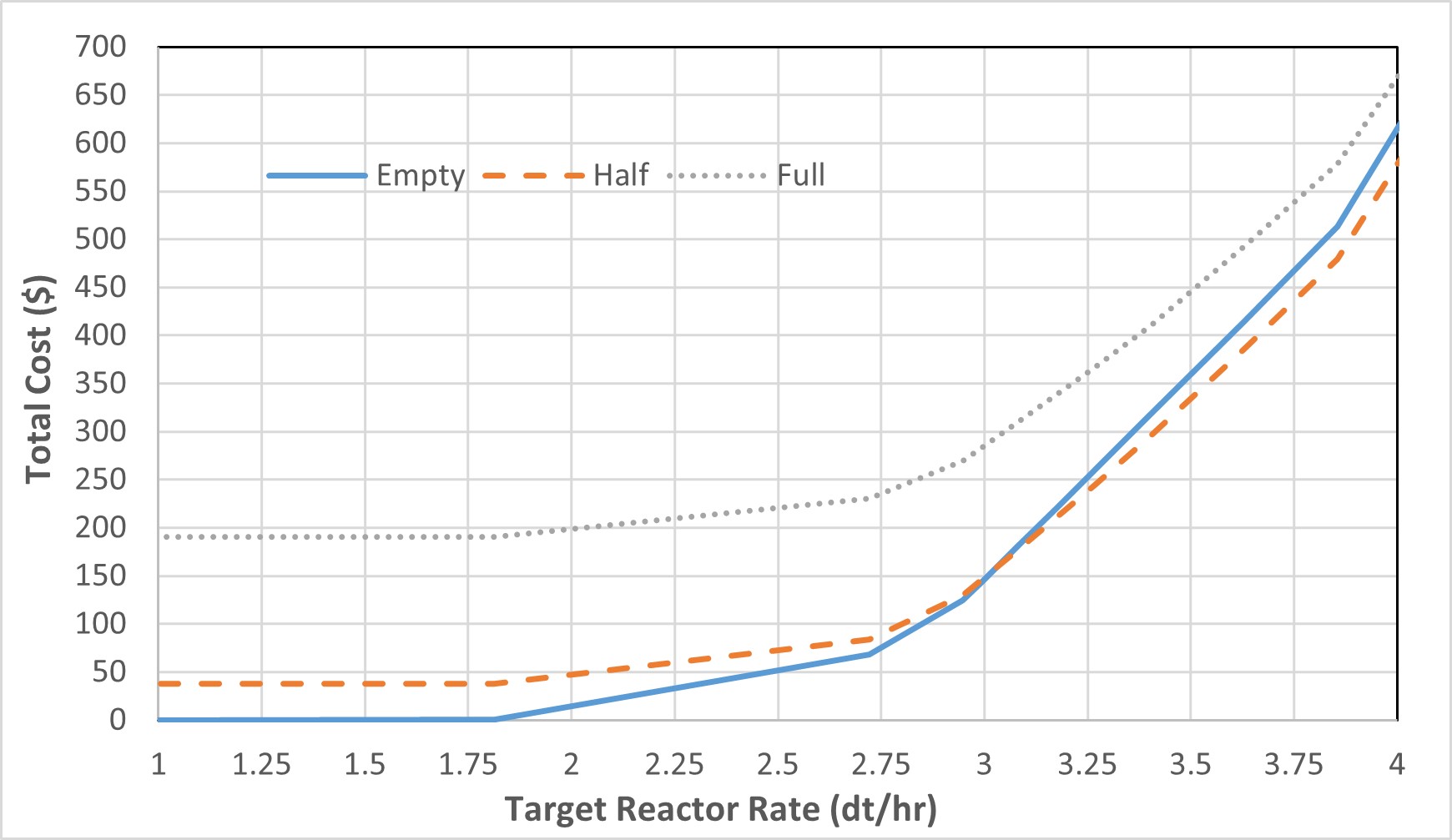}
    	\caption{Sensitivity analysis for initial inventory level.}
    	\label{fig:sensitivity_Initial_Inv}
    \end{subfigure} 
    \hfill
    \begin{subfigure}[t]{0.48\textwidth}
        \centering
    	\includegraphics[width=\textwidth]{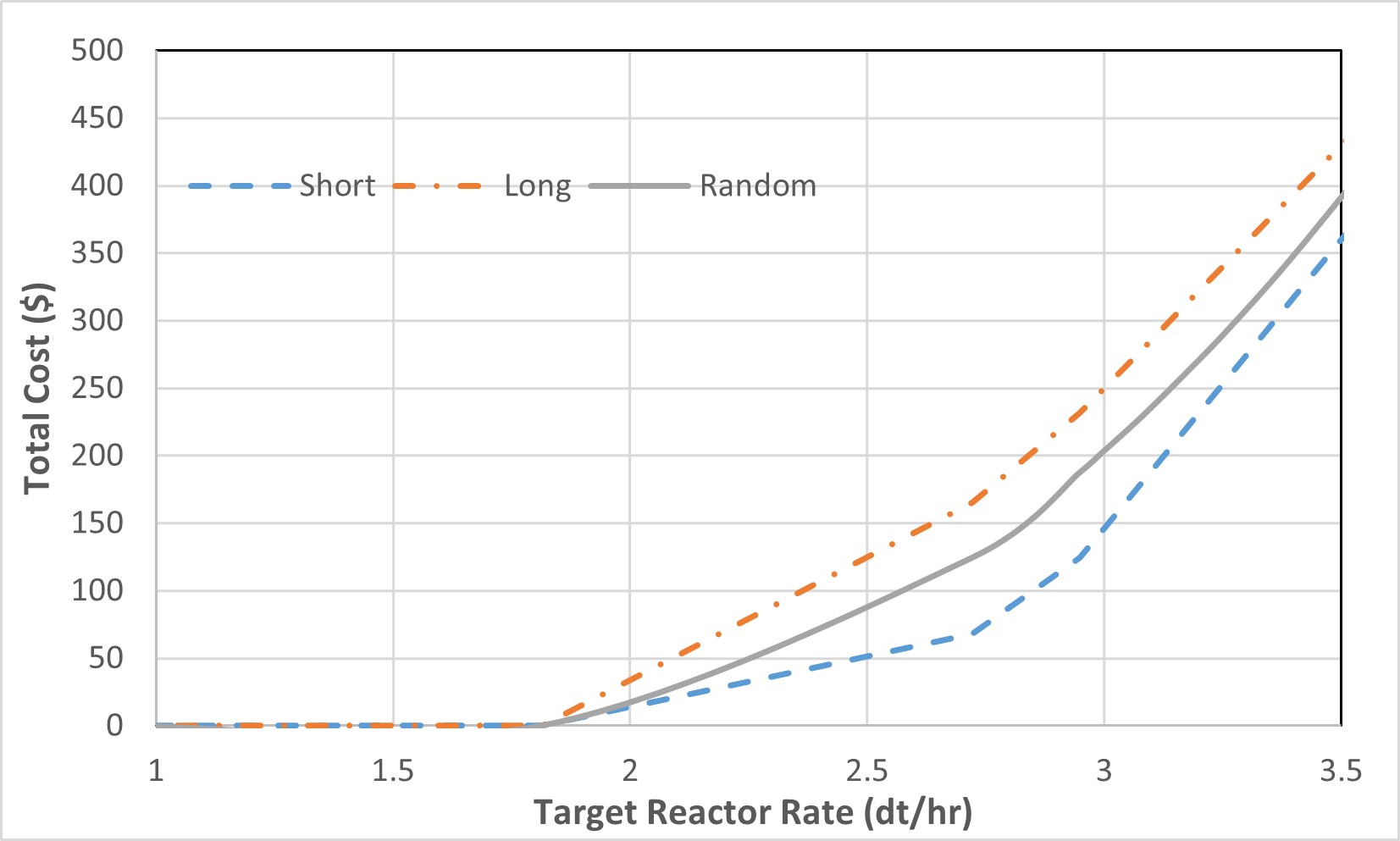}
    	\caption{Sensitivity analysis for the length of the bale sequencing pattern.}
    	\label{fig:sensitivity_SequenceL}
    \end{subfigure}
    \caption{Sensitivity analysis}
\end{figure}

Figure~\ref{fig:sensitivity_Initial_Inv} shows how the average operating cost changes as the target reactor rate increases from $1.0\ dt/hr$ to $4.0\ dt/hr$ for the three initial inventory levels considered in our sensitivity analysis, using the base-case problem described in Section~\ref{sec:model-evaluation}. Note that the total cost calculation only considers the total inventory holding cost accumulated over the planning horizon but does not include the cost for the initial inventory.

In Figure~\ref{fig:sensitivity_Initial_Inv}, we observe that for smaller target rates, starting the operation with empty storage equipment is more economical. Since the system can achieve the target reactor rate without additional inventory and the system incurs inventory holding cost for not used inventory. However, as the target reactor feeding rate exceeds $2.95-3.00\ dt/hr$, having a half-full initial inventory becomes more economical. The additional initial inventory makes it easier to achieve the target reactor feeding rate. We also see that operating with a full initial inventory is always more costly than the other two strategies for the range of target reactor rates considered. 

Another interesting observation from these experiments was that, when operating with a half or full initial inventory, the static solution obtained from the deterministic mean-value problem is infeasible. In order to satisfy the inventory levels set by this static solution, we found out that the system needs to feed the pelleting mill with a rate higher than the pelleting mill's infeed capacity. This observation points to the need for developing a SP model. 

\subsection{Sensitivity Analysis: Length of the Bale Sequencing Pattern} \label{sec:Sensitivity_SequenceLen}
In this section, we analyze the impact of the sequencing choice on the operational cost and reactor utilization. For this analysis, we introduce two different sequencing strategies: the long sequence and the random sequence. The long sequence considers that the system processes every bale of a particular moisture level before processing bales of another moisture level. The operator forms the random sequence by picking bales randomly. In our experiments, we considered 10 random sequences and presented the average results. We solve the base-case problem described in Section~\ref{sec:model-evaluation} for the three different sequencing strategies (short, long, and random) with a range of target reactor rates.

Figure~\ref{fig:sensitivity_SequenceL} clearly shows the advantage of the short sequence over the entire range of target reactor rates, and in certain cases, it results in up to a $59\%$ reduction in the total cost compared to alternative options. This is because the repeating pattern in the short sequence enables the system to accumulate inventory while processing low-moisture bales and to consume the inventory while processing high-moisture bales. The performance of the random sequence lies in between the short and long sequencing strategies. One caveat is that our model does not consider the cost of sequencing the bales. Therefore, if the cost of creating a short sequence (due to additional operations such as sorting) is high, a biomass processing plant may simply use the random sequence as a rule of thumb in their operation.

\subsection{Sensitivity Analysis: Sequence Order} \label{sec:Sensitivity_SequenceOrder}
In this section, we analyze the impact of the sequencing order on the total cost and reactor utilization. In the base-case problem, we use the short sequence that repeats the pattern of processing one high-moisture bale, one medium-moisture bale, and three low-moisture bales. Let us call this ordering the high-start sequence. Next, we will compare this the low-start sequence, which we defined to have the repeating pattern of three low-moisture bales, one medium-moisture bale, and one high-moisture bale. We use the base-case problem described in Section~\ref{sec:model-evaluation} with the high-start and low-start sequences for a range of target rates for the reactor.

\begin{figure}[h]
    \centering
    \begin{subfigure}[t]{0.48\textwidth}
        \centering
    	\includegraphics[width=\textwidth]{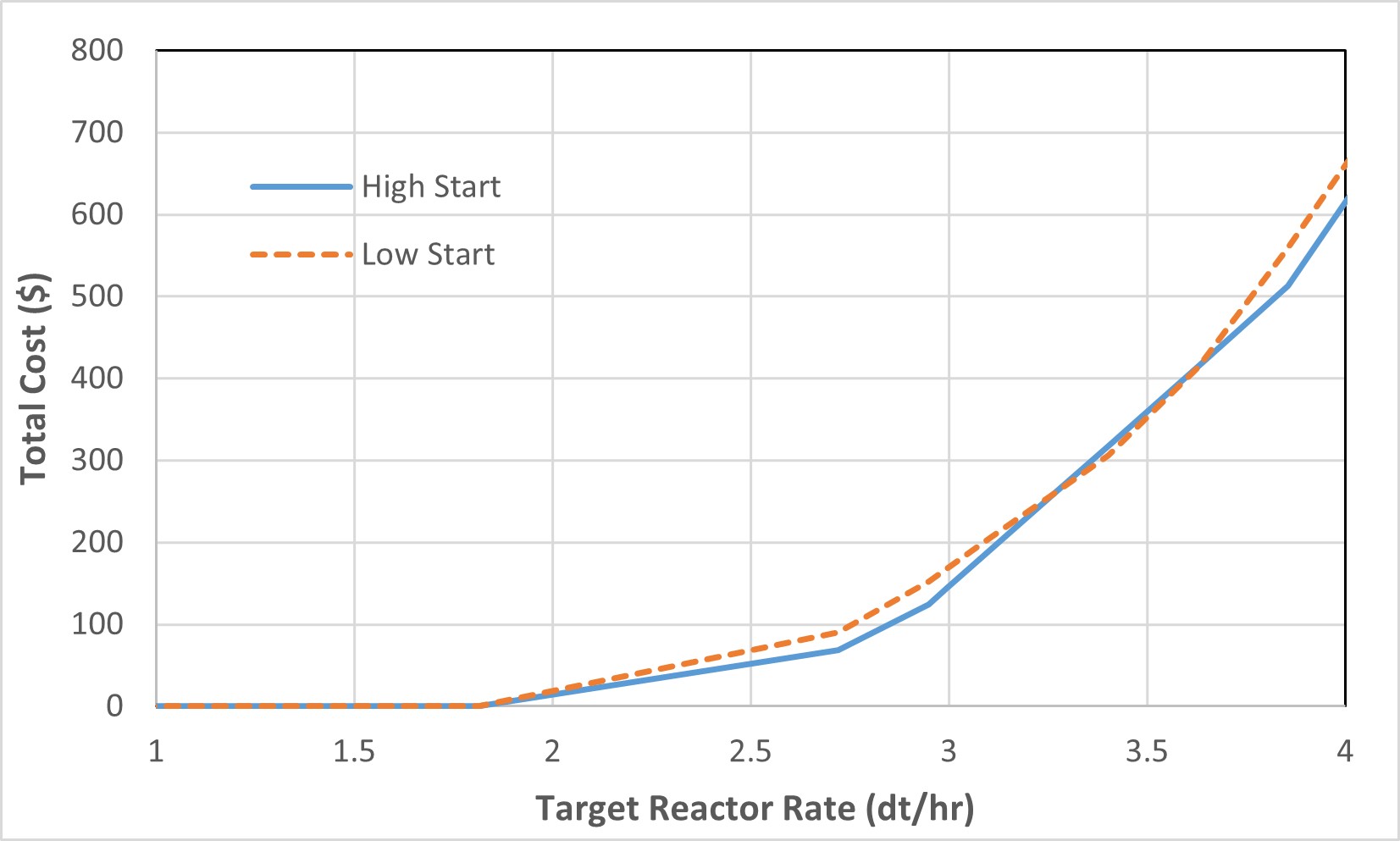}
    	\caption{Total Cost for different target rates.}
    	\label{fig:sensitivity_SequenceOrder}
    \end{subfigure} 
    \hfill
    \begin{subfigure}[t]{0.48\textwidth}
        \centering
    	\includegraphics[width=0.9\textwidth]{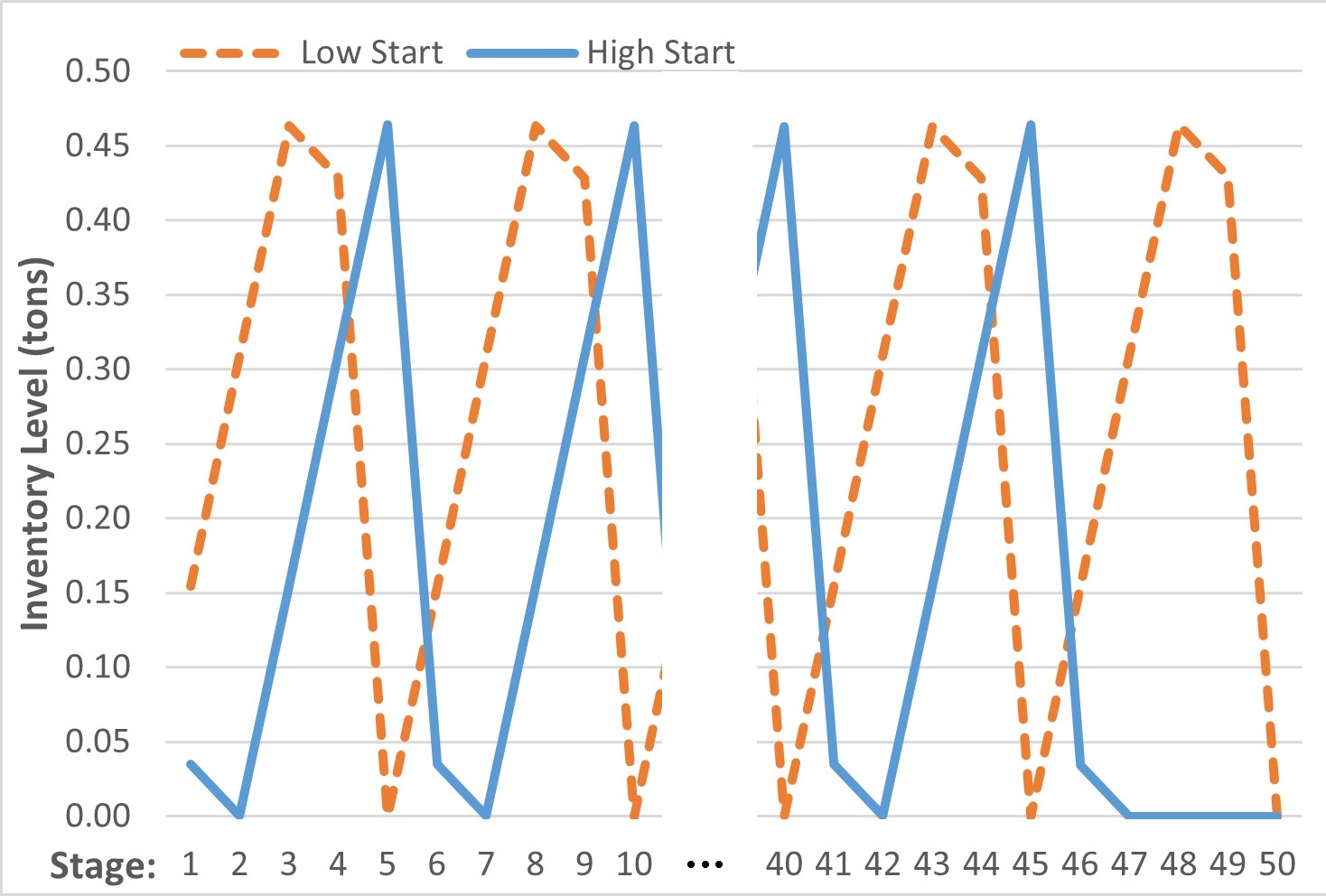}
    	\caption{Inventory levels for $r\ =\ 2.95\ dt/hr$.}
    	\label{fig:sensitivity_Inv_SequenceOrder}
    \end{subfigure}
    \caption{Sensitivity analysis for the bale sequence order.}
\end{figure}

We observe from Figure~\ref{fig:sensitivity_SequenceOrder} that for lower target rates, starting with high-moisture bales results in a better performance. To help facilitate the understanding of this phenomenon, we show in Figure~\ref{fig:sensitivity_Inv_SequenceOrder} the inventory levels for high-start and low-start strategies in each stage. We see that the inventory level presents a cyclic pattern for both strategies. At the beginning and at the end of the planning horizon, the low-start strategy accumulates inventory that will be used when bales with high moisture level are processed. This results in a $56\%$ higher inventory level on average. As the target rate increases, the low-start strategy performs better.

\subsection{Key Findings}\label{sec:Sensitivity_summary}
This section summarizes the take-away messages from our numerical experiments and sensitivity analysis. First, using multi-stage SP to model sequential decisions of biorefinery operations improves reactor utilization and costs compared to alternative models (two-stage SP, deterministic mean-value problem) that use only static solutions. Second, the bale-by-bale definition of stages appears to yield the best performance.
A more detailed stage definition results in only a marginal gain at the cost of a 217.7\% increase in computational time. The combined stage definition results in a 28.4\% reduction in computational time but at the expense of a 33.5\% deterioration in the solution quality. Third, filling the storage equipment to its half capacity as an initial inventory level improves the system's capability to meet the reactor's target feeding rate and reduces costs for high reactor processing rates. Finally, we observe that sequencing bales based on their moisture levels can lead to as much as 59\% cost reduction for certain target reactor rates. The short sequences of low-, medium-, and high-moisture bales performs best. 

\section{Summary and Conclusion}\label{sec:Conclusion}
The main objective of this work is to minimize operational costs and improve the reliability and flexibility of a biomass processing system in a biorefinery to changing conditions of the processed biomass. We propose a multi-stage SP model that identifies the initial inventory level, initial equipment settings and creates a policy for future operational decisions. The uncertain moisture content impacts the processing speed of equipment, biomass density and the volumetric capacity of storage space. We take a system's approach to this problem and use a network flow model to characterize the flow of biomass in the system. We evaluated the system performance based on inventory holding costs and achieving target reactor rate. We ensure the consistent high utilization of the reactor by penalizing the violation of the target reactor feeding rate. 

This research makes a number of contributions to the literature. First, we justify the modeling choice and present the value of multi-stage SP to optimize biomass processing operations under uncertainty. Second, we develop the proposed model using real-life data about biomass pre-processing. The model enables the integration of sensor technologies into automated system control. As a result, we expect that practitioners will find the results of this study useful and inspiring for the transitioning into Industry 4.0. Third, we make several important observations. We observe that revising the control decisions every time the system starts processing a new bale provides a significant cost advantage compared to less frequent control schemes. We observe that when sequencing of biomass bales is cheap, practitioners should sequence the bales based on moisture content in patterns that repeat frequently. If sequencing is costly, randomly picking biomass bales may be a good strategy.

We have identified several research directions that are worth exploring in the future. First, the model can be extended to consider other uncertain biomass characteristics such as particle size, carbohydrate, and ash contents. Considering these biomass characteristics will increase the variation among scenarios and highlight the impacts of sequencing and inventory holding strategies. In addition, the inclusion of carbohydrate and ash content might require the identification of critical constraints regarding biomass quality. Second, the proposed model can be extended to consider multiple biomass types and include the sequencing strategies as decision variables. This will lead to endogenous uncertainty as the biomass characteristics during the process depend on the sequencing decision. This more challenging problem will motivate the development of new approximation schemes and computationally efficient solution algorithms.

\bibliographystyle{tfcad}
\bibliography{references2}

\newpage

\section{Appendices}
\appendix
\section{Data Tables}\label{App:DataTables}

\begin{table}[H]
\centering
	\scriptsize
	{%
        	\begin{tabular}{lcccc}
    			\toprule
    			& \multicolumn{2}{c}{\textbf{Biomass Processed}}& \multicolumn{2}{c}{\textbf{Biomass Processed}}\\
    			& \multicolumn{2}{c}{\textbf{in Grinder 1}}& \multicolumn{2}{c}{\textbf{in Grinder 2}}\\
    			\cmidrule(lr){2-3}
    			\cmidrule(lr){4-5}
    			\textbf{Moisture} & $\boldsymbol{\rho^{50}}$& \multirow{2}{*}{$\boldsymbol{\rho^{90}/\rho^{10}}$} & $\boldsymbol{\rho^{50}}$& \multirow{2}{*}{$\boldsymbol{\rho^{90}/\rho^{10}}$} \\
    			\textbf{ \hspace{.1mm} Level}&$(mm)$ & &$(mm)$ &\\
    			\cmidrule(lr){1-1}
    			\cmidrule(lr){2-2}
    			\cmidrule(lr){3-3}
    			\cmidrule(lr){4-4}
    			\cmidrule(lr){5-5}
    			Low		&  $ 1.95 $ & $ 12.5 $ & $ 0.65 $& $ 6.50 $ \\
    			Medium	&  $ 2.35$ & $ 12.0$ & $ 0.70$& $ 7.50$ \\
    			High	&  $ 1.75$ & $ 10.0$ & $ 0.60 $& $ 9.50 $ \\			
    			\bottomrule
    	\end{tabular}}
    	\captionof{table}{Particle size distribution percentiles.}
    	\label{tab:PSD}
\end{table}%

\begin{table}[H]
\centering
	\scriptsize
	{%
        	\begin{tabular}{cccc}
    			\toprule		
    			\textbf{Moisture}& \multirow{2}{*}{\textbf{Low}}& \multirow{2}{*}{\textbf{Med}}& \multirow{2}{*}{\textbf{High}}\\
    			\textbf{Level}& & & \\
    			\cmidrule(lr){1-1}
    			\cmidrule(lr){2-2}
    			\cmidrule(lr){3-3}
    			\cmidrule(lr){4-4} 
    			\textbf{Bypass} & \multirow{2}{*}{85.7}&\multirow{2}{*}{81.1} &\multirow{2}{*}{92.8} \\
    			\textbf{Ratio $(\%)$} & & & \\
    			\bottomrule
    	\end{tabular}}
    	\captionof{table}{Secondary grinder bypass ratios.}
    	\label{tab:Bypass}
\end{table}%

\begin{table}[H]
\centering
	\scriptsize
	{%
        	\begin{tabular}{lccc}
        			\toprule		
        			& \multicolumn{3}{c}{\textbf{Moisture Level}}\\
        			\cmidrule(lr){2-4}
        			& \textbf{Low}& \textbf{Med}& \textbf{High}\\
        			\textbf{Equipment}& $(dt/ hr)$& $(dt/ hr)$& $(dt/ hr)$\\
        			\cmidrule(lr){1-1}
        			\cmidrule(lr){2-2}
        			\cmidrule(lr){3-3}
        			\cmidrule(lr){4-4}
        			Grinder 1 & 5.23&4.53 &2.20 \\
        			Grinder 2 &5.23 &2.80 &1.59 \\
        			Pelleting &4.76 & 3.81& 3.34\\
        			Reactor & 4.81 & 4.81 & 4.81\\
        			\bottomrule
        			\\
        	\end{tabular}}
        	\captionof{table}{Equipment infeed rate limits.}
        	\label{tab:infeed_limit}
\end{table}%

\begin{table}[H]
	\scriptsize
	\centering
	{%
	    \setlength{\tabcolsep}{0.6em}
		\begin{tabular}{lccccc}
			\toprule
			
			& \multicolumn{3}{c}{\textbf{Moisture Loss}}& \multicolumn{2}{c}{\textbf{Dry Matter Loss}}\\
			\cmidrule(lr){2-4}
			\cmidrule(lr){5-6}
			\textbf{Moisture} & \textbf{Grinder 1}& \textbf{Grinder 2} & \textbf{Pelleting}& \textbf{Grinder 1} & \textbf{Grinder 2}\\
			\textbf{ \hspace{.1mm} Level}&$(\%)$ &$(\%)$ &$(\%)$ &$(\%)$ &$(\%)$\\
			\cmidrule(lr){1-1}
			\cmidrule(lr){2-2}
			\cmidrule(lr){3-3}
			\cmidrule(lr){4-4}
			\cmidrule(lr){5-5}
			\cmidrule(lr){6-6}
			Low		&  $0.50$ & $0.7 0$ & $0.00$ & $1.50$ & $0.50$ \\
			Medium	&  $3.00$ & $3.00$ & $1.50$ & $1.50$ & $0.50$ \\
			High	&  $4.77$ & $4.00$ & $3.90$ & $1.50$ & $0.50$ \\			
			\bottomrule
	\end{tabular}}
	\caption{Moisture and dry matter changes.}
	\label{tab:processing_changes}
\end{table}%

\section{Regression Analysis for Biomass Density}
\label{sec:Regress_App}

\begin{table}[H]
	\scriptsize
	\centering
	{%
	    \setlength{\tabcolsep}{0.6em}
		\begin{tabular}{lcccc}
			\toprule			
			& & \multicolumn{3}{c}{\textbf{P-value}}\\
			\cmidrule(lr){3-5}
			\textbf{Regression} & $\mathbf{R^2}$ & $\mathbf{\alpha_i^0}$ &$\mathbf{m_{it}}$ & $\mathbf{\rho^{50}_i}$\\
			\cmidrule(lr){1-1}
			\cmidrule(lr){2-2}
			\cmidrule(lr){3-3}
			\cmidrule(lr){4-4}
			\cmidrule(lr){5-5}
			(16) & 0.956 & $4.1*10^{-142}$ & $9.5*10^{-134}$ & $1.3*10^{-69}$ \\
			(17) & 0.945 & $2.1*10^{-69}$ & $2.4*10^{-75}$ & $2.5*10^{-39}$\\
			\bottomrule
	\end{tabular}}
	\caption{Regression Analysis Statistics}
	\label{tab:reg_stats}
\end{table}%

\end{document}